\documentstyle [twocolumn,amssymbols]{article}

\renewcommand{\theequation}{\thesection.\arabic{equation}}

\newcommand{\for}[1]{\mbox{(\ref{#1})}}

\newtheorem{defn}{Definition}[section]
\newtheorem{lem}[defn]{Lemma}
\newtheorem{prp}[defn]{Proposition}
\newtheorem{thm}[defn]{Theorem}
\newtheorem{cor}[defn]{Corollary}

\newenvironment{dfn}{\begin{defn}\rm }{\end{defn}}
\newcommand{\epr}{\begin{flushright} $\Box$ \end{flushright}}
\newenvironment{prf}{{\em Proof}\newline }{\epr}

\newcommand{\note}{{\em Note. }}
\newcommand{\notes}{{\em Notes }}
\newcommand{\remark}{{\em Remark. }}
\newcommand{\remarks}{{\em Remarks }}

\newcommand{\pdo}{PDO}
\newcommand{\ado}{A$\Delta $O}
\newcommand{\pdos}{PDO's}
\newcommand{\ados}{A$\Delta $O's}
\newcommand{\ra}{\rightarrow}
\newcommand{\eps}{\varepsilon}
\newcommand{\pa}{\partial}

\newcommand{\ch}{{\rm ch}}

\begin{document}
\rightline{June 1993}
\rightline{funct-an/9306002}

\vspace{-5.5ex}
{\noindent\Large\bf Commuting Difference Operators \\
  with Polynomial Eigenfunctions}

\vspace{3ex}
{\noindent\large J.F. van Diejen}\footnote{E-mail
address: jand@fwi.uva.nl}

{\noindent\small
Department of Mathematics and Computer Science, University of Amsterdam,\\
Plantage Muidergracht 24, 1018 TV Amsterdam, The Netherlands.}

\vspace{1ex}
\small\noindent {\bf Abstract.}
We present explicit generators $\hat{D}_1,\ldots ,\hat{D}_n$
of an algebra of commuting difference operators
in $n$ variables with trigonometric coefficients.
The algebra depends, apart from two scale factors,
on five parameters.
The operators are simultaneously diagonalized by Koornwinder's
multivariable generalization of the Askey-Wilson polynomials.
For special values of the parameters and via limit transitions, one obtains
difference operators for the Macdonald polynomials that are associated with
(admissible pairs of) the classical root systems:
$A_{n-1}$, $B_n$, $C_n$, $D_n$ and $BC_n$.
By sending the step size of the differences to zero, the difference operators
reduce to known hypergeometric differential operators.
This limit corresponds to sending $q\ra 1$; the eigenfunctions reduce to
the multivariable Jacobi polynomials of Heckman and Opdam.
Physically the algebra can be interpreted as
an integrable quantum system that generalizes the (trigonometric)
Calogero-Moser systems related to classical root systems.

\vspace{1ex}
{\sl Key words \& Phrases:} commuting difference operators, (joint)
eigenfunctions, multivariable q-polynomials, classical root systems,
integrable quantum systems.

\tableofcontents

\setcounter{equation}{0}
\addtocontents{toc}{\protect\vspace{-4ex}}
\markright{1 Introduction}
\section{Introduction}\label{sect1}
Over the past few years, progress has been made with the study of orthogonal
polynomials in more than one variable.
It has turned out that a lot of classical results concerning orthogonal
polynomials depending on only one variable admit generalization to many
variables.
Such generalizations can be viewed naturally in a Lie-theoretic setting:
for each root system $R$, there exist associated families of multivariable
polynomials. The number of variables coincides with the rank $n$ of the root
system. Families of multivariable Jacobi polynomials related to root systems
have been studied by Heckman and Opdam \cite{heck1,heck2}.
Recently, a more elementary account of some of these results was presented
in \cite{heck3}.
For $R=BC_1$ the Heckman-Opdam-Jacobi polynomials reduce to the classical
Jacobi polynomials in
one variable.
In a yet unpublished manuscript, Macdonald has introduced q-versions of the
Heckman-Opdam families \cite{mac1}.
See \cite{mac2,koo3} for a summary of these results and \cite{mac1.5} for
lectures devoted to the special case $R=A_n$.
If $R=BC_1$, then Macdonald's polynomials coincide with
the continuous q-Jacobi polynomials. (For information
on continuous q-Jacobi polynomials see e.g. \cite{akw}).
Recently, Macdonald's results pertaining to the root system $BC_n$ have been
generalized by Koornwinder \cite{koo3}. He finds $BC_n$-type multivariable
versions of the Askey-Wilson polynomials \cite{akw}. Again, the one-variable
case is recovered by specializing to $BC_1$.

A crucial ingredient in the construction of the above families is the
existence of an operator of which the polynomials are eigenfunctions.
In case of the Heckman-Opdam-Jacobi polynomials this operator is a
second order partial
differential operator (\pdo{}) which is named hypergeometric differential
operator.
For Macdonald's and Koornwinder's polynomials the relevant operator is an
analytic difference operator (\ado{}).
A similarity transformation turns the hypergeometric \pdo{} into
a differential operator which is self-adjoint with respect to Lebesgue
measure.
{}From the perspective of physics, one can look upon the latter \pdo{}
as being the Hamiltonian of a quantum system of
$n$ particles in dimension one. Such quantum systems have been studied for
quite some time in the physics literature; they are known as
(generalized) Calogero-Moser systems \cite{cal,sut1,sut2,op2}.
The systems studied in \cite{cal} and \cite{sut1,sut2} correspond to the root
system $A_{n-1}$.  See the survey paper \cite{op2} for the generalization to
arbitrary root systems.
The motion of the corresponding classical systems has also been
studied \cite{mos} ($R=A_{n-1}$) and \cite{op1} (arbitrary $R$).

It has been shown by Heckman and Opdam that the hypergeometric differential
operator is but
one member of an algebra of commuting \pdos{} having
the multivariable Jacobi
polynomials as their joint eigenfunctions \cite{heck1,heck2,opd1,opd2}.
This algebra is generated by $n$ independent \pdos{}.
This state of affairs can be expressed by saying that
the corresponding (generalized)
Calogero-Moser system is quantum integrable. For $R=A_{n-1}$,
Liouville integrability of the classical system, i.e. the existence of
$n$ independent integrals in involution,
was already proved by Moser \cite{mos} using a Lax pair formulation.
For a partial generalization of
this result to the other classical (i.e. non-exceptional) root systems,
see \cite{op1} and \cite{ino}.

Just as the hypergeometric \pdo{}, Macdonald's difference
operator for $R=A_{n-1}$ is related to certain known
quantum systems of $n$ particles.
The systems of interest were originally introduced as a relativistic
generalization of the Calogero-Moser systems \cite{rui1} (classical) and
\cite{rui2}
(quantum). (See \cite{rui3} for a survey and connections with certain
soliton PDE's and exactly solvable quantum field theories).
The relativistic generalization of the Calogero-Moser system is (quantum)
integrable too:
explicit formulas representing $n$ independent commuting
integrals are presented in \cite{rui1,rui2}.
At the quantum level these integrals are \ados{}, which are related
to Macdonald's $A_{n-1}$-type difference
operators $E_{\omega_r}$ via a similarity transformation \cite{koo0}.
The operators $E_{\omega_r}$
(associated with the fundamental weights $\omega_r$)
generate an algebra of commuting \ados{} which have
Macdonald's $A_{n-1}$-type polynomials as joint eigenfunctions.

Also for root systems other than $A_{n-1}$ one expects that there exist
algebras of commuting
difference operators that are simultaneously diagonalized by
the Macdonald polynomials.
The purpose of the present paper is to introduce such difference
operators for the Koornwinder polynomials.
We will show that via limit transitions and/or specialization of parameters,
this also leads to the corresponding
\ados{} for those families of Macdonald polynomials that are connected with
(admissible pairs of) the classical series: $A_{n-1}$, $B_n$, $C_n$, $D_n$ and
$BC_n$;
the Koornwinder polynomials then reduce to the latter
Macdonald polynomials.
By sending the step size of the differences to zero (this corresponds to the
limit $q\ra 1$) our \ados{} go over in \pdos{}. Thus, we recover the
commuting hypergeometric \pdos{} associated with the classical
root systems as a limit case.

Our difference operators constitute a new integrable quantum system
of $n$ particles in dimension one.
In this paper, however, we will not pay much attention to
this interpretation of the \ados{}; instead we will
emphasize the connection with orthogonal polynomials.
In a forthcoming paper the author intends to return to the
question of integrability of these and related $n$-particle systems,
at the level of both quantum and classical mechanics \cite{die}.
In particular, possible generalization to integrable systems consisting of
commuting difference operators
with elliptic functions as coefficients will be discussed in \cite{die}.

Before outlining the contents of this paper in more detail, let us
mention two more connections of interest.
For special values of the parameters (namely those corresponding
to root multiplicities), the
system of hypergeometric \pdos{} coincides with the radial reduction of the
algebra of invariant differential operators on certain symmetric spaces $G/K$
\cite{op2,heck1}. It seems natural to ask oneself the question whether, for
special
values of the parameters, our system of difference operators can be seen in
some
way as radial reduction of certain \ados{} connected with quantum homogeneous
spaces. Recent results on the quantum group interpretation of Macdonald's
$A_{n-1}$-type polynomials \cite{koo2,nou} indeed seem to point in this
direction. However, apart from this special case no relations of this
kind are known to the author.

Recently, Cherednik introduced commuting difference operators connected
with Knizh\-nik-Zamolodchikov-type difference equations
associated with affine root systems \cite{ch1,ch2}.
It is claimed by the author that for $A_{n-1}$ his operators
coincide with the Macdonald \ados{}. It would be interesting to investigate
the relation with the explicit \ados{} of the present paper.
Specifically, one would like to know whether the Cherednik operators that
are associated with the classical root systems are all special cases of
the latter \ados{}.

The paper is organized as follows:
in Section~\ref{sect2} we introduce $n$ independent \ado{} 's
$\hat{D}_r$, $r=1,\ldots ,n$ in the (real) variables $x_1, \ldots ,x_n$.
These \ados{} (which we write down explicitly) depend,
apart from two scale factors, on five parameters.
The simplest operator, viz. $\hat{D}_1$, coincides
(up to an irrelevant multiplicative constant) with Koornwinder's difference
operator $D_{\eps_1}$; for $r\geq 2$ our operators are new. By setting a
certain parameter to zero, $\hat{D}_r$ reduces to the $r$th elementary
symmetric function of the operators $\hat{D}_1(x_j)$, $j=1,\ldots ,n$.
(Here $\hat{D}_1(x_j)$ denotes
the one-variable version of $\hat{D}_1$ with $x_j$ as variable).

Section~\ref{sect3} contains the main results of the paper.
In subsequent subsections it is shown that:
{\em i.} $\hat{D}_r$ leaves invariant certain finite-dimensional highest weight
spaces (triangularity);
{\em ii.} $\hat{D}_r$ is symmetric with respect to the $L^2$ inner
product with weight function $\Delta$ ($\Delta dx$ being the orthogonality
measure of Koornwinder's polynomials).
Combination of these two facts implies that
$\hat{D}_r$ is diagonalized by the Koornwinder polynomials.
We would like to mention that this method to diagonalize the difference
operators $\hat{D}_r$ resembles very much the approach that was originally used
by Sutherland, which leads to the spectrum and eigenfunctions of the
(trigonometric) Calogero-Moser system \cite{sut1,sut2}.
The spectrum of the
operators is computed explicitly. As a result, we obtain a Harish-Chandra-type
isomorphism between the abelian algebra generated by
$\hat{D}_1,\ldots ,\hat{D}_n$
and the symmetric algebra in $n$ variables.
For $\hat{D}_1$ the discussion in this section amounts to a reproduction
of results already obtained by Koornwinder \cite{koo3}.
The difference between our presentation (when restricted to the case $r=1$)
and that of \cite{koo3} is that (by exploiting a calculation of residues
and some asymptotics) we avoid certain (rather long) calculations in
Koornwinder's paper leading to the triangularity and the eigenvalues of
$\hat{D}_1$.

In Section~\ref{sect4} we discuss
the transition to Heckman and Opdam's commuting hypergeometric \pdos{}
related to the root system $BC_n$.
In this limit, which amounts to sending $q\ra 1$, the
eigenfunctions converge to the $BC_n$-type Jacobi polynomials.

In Section~\ref{sect5} we study various special cases related to
the classical root systems. In 5.2 we introduce a limit transition
leading to the $A_{n-1}$ root system, which is the most interesting
from a physical viewpoint. Specifically, we show how the $A_{n-1}$
Macdonald polynomials and \ados{} can be obtained from the Koornwinder
polynomials and the \ados{} $\hat{D}_1,\ldots ,\hat{D}_n$, respectively.
This novel transition can be applied to other situations as well.
(For example, it enables one to view the $A_{n-1}$
Jacobi polynomials as limits of their $D_n$ counterparts).

In 5.3 and 5.4 we show how the Macdonald polynomials associated
with the remaining root systems can be obtained from the Koornwinder
polynomials by suitable specialization of the parameters. In contrast to
our account for the $A_{n-1}$ case, which is quite self-contained, this
involves
various concepts from \cite{mac1}. We have attempted to render our account more
accessible by collecting some preliminaries in 5.1; this subsection
can be skipped at first reading and referred back to as needed.

\setcounter{equation}{0}
\addtocontents{toc}{\protect\vspace{-2ex}}
\markright{2 Introducing the Difference Operators}
\section{Introducing the Difference Operators}\label{sect2}
In this section the operators $\hat{D}_r$, $r=1,\ldots , n$, are introduced and
their combinatorial structure is discussed.

\subsection{The Operator $\hat{D}_r$}\label{ados}
In order to write down our difference operators we first introduce some
notation. Let $v_a(z)$ and $v_b(z)$ be the following trigonometric
functions:
\begin{eqnarray}
v_a(z) &=& \frac{\sin\alpha (\mu +z)}{\sin\alpha z}, \label{va} \\
v_b(z) &=& \frac{\sin\alpha (\mu_0 +z)}{\sin\alpha z}
           \frac{\cos\alpha (\mu_1 +z)}{\cos\alpha z}
   \frac{\sin\alpha (\mu_0^\prime +\gamma +z)}{\sin\alpha (\gamma +z)}
   \frac{\cos\alpha (\mu_1^\prime +\gamma +z)}{\cos\alpha (\gamma +z)},
\label{vb}
\end{eqnarray}
with $\alpha$, $\gamma$ and $\mu$, $\mu_\delta$, $\mu_\delta^\prime$,
($\delta =0,1$) complex parameters. For later purposes, it is convenient
to parametrize $\gamma$ according to:
\begin{equation}
\gamma \equiv i\beta /2.\label{gambet}
\end{equation}
We form the following multivariable functions using $v_a$ and $v_b$ as
elementary constituents:
\begin{eqnarray}
V_{\eps J;K} &\equiv&
\prod_{j\in J} v_b(\eps_j x_j)
\prod_{\stackrel{j,j^\prime \in J}{j<j^\prime}}
v_a(\eps_j x_j+\eps_{j^\prime} x_{j^\prime})\,
v_a(\eps_j x_j+\eps_{j^\prime} x_{j^\prime} +2\gamma ) \nonumber\\
& &\times
\prod_{\stackrel{j\in J}{k\in K}} v_a(\eps_j x_j +x_k)\, v_a(\eps_j x_j -x_k),
\label{VJ}
\end{eqnarray}
with
\begin{equation}
J,K \subset \{ 1,\ldots , n\} ,\;\;\;\;
J\cap K =\emptyset;\;\;\;\;\;\;\;
\eps_j \in \{ +1,-1\} .
\end{equation}
The variables $x_1,\ldots , x_n$ are assumed to be real.
The function $V_{\eps J;K}$ depends on the index sets $J$, $K$ and on a
collection of prescribed signs $\eps_j$, $j\in J$; it serves as a
building block from which the coefficients of $\hat{D}_{r}$ are constructed.
The \ados{} read explicitly
\begin{eqnarray}
\hat{D}_{r} \equiv
\sum_{\stackrel{J\subset \{ 1,\ldots , n\} ,\, |J|=r}{\eps_j=\pm 1,\, j\in
J}}\;\;\;
\sum_{\stackrel{\emptyset\subsetneqq J_1\subsetneqq \cdots \subsetneqq J_s=J}
               {1\leq s\leq r} }
(-1)^{s-1}\prod_{1\leq s^\prime \leq s}
V_{\eps (J_{s^\prime}\setminus J_{s^\prime -1});J_{s^\prime}^c}
\left( e^{-\beta \hat{\theta}_{\eps J_1 } }-1 \right) ,\label{Dr1} \\
\hfill r=1,\ldots , n, \nonumber
\end{eqnarray}
with $J_0 \equiv \emptyset$ and
\begin{equation}
\hat{\theta}_{\eps J}\equiv \sum_{j\in J}\eps_j \hat{\theta}_j
,\;\;\;\;\;\;\;\;
\hat{\theta}_j\equiv \frac{1}{i}\frac{\pa}{\pa x_j}.\label{shift}
\end{equation}

\vspace{1ex}
\noindent\remarks
{\em i.} $|J|$ denotes the cardinality of $J$ and $J^c$ is the complement set
$\{ 1,\ldots , n\} \setminus J$. \\
{\em ii.} The first summation in Eq. \for{Dr1} is over all index sets
$J\subset \{ 1,\ldots , n\} $ with cardinality $r$ and over all flippings of
the signs
$\eps_j \in \{ +1,-1\}$, $j\in J$; the second summation is over all strictly
increasing sequences of subsets in $J$:
\begin{equation}
\emptyset\subsetneqq J_1\subsetneqq J_2\subsetneqq\cdots
\subsetneqq J_{s-1}\subsetneqq
J_s=J,\;\;\;\;\;\;  1\leq s \leq |J|. \label{part}
\end{equation}
{\em iii.} The exponential $\exp (-\beta \hat{\theta}_j)$ acts on a
function $f(x_1,\ldots , x_n)$, which is analytic in the variables
$x_1,\ldots , x_n$, as a (complex) shift:
\begin{equation}
(e^{-\beta \hat{\theta}_j}f)(x_1,\ldots , x_n)=
f(x_1,\ldots , x_{j-1},x_{j}+i\beta ,x_{j+1},\ldots , x_n).
\end{equation}
Hence, $\hat{D}_r$ is indeed an analytic difference
operator. \\
{\em iv.} The `hats' in Eqs. \for{Dr1}, \for{shift} are used to emphasize that
$\hat{D}_r$ and $\hat{\theta}_j$ are operators rather than ordinary (complex)
functions or variables.\\
{\em v.} In the simplest case, i.e. for $r=1$, Eq. \for{Dr1} reduces to
\begin{equation}
\hat{D}_1 = \sum_{\stackrel{1\leq j\leq n}{\eps =\pm 1}}
v_{b}(\eps x_j)\prod_{k\neq j} v_a (\eps x_j+x_k) v_a (\eps x_j-x_k)
\left( e^{-\eps\beta \hat{\theta}_j}-1 \right) .\label{D1}
\end{equation}
{\em vi.} It is clear that the operator $\hat{D}_r$ is invariant under
permutations of
the variables $x_1,\ldots , x_n$. Furthermore, the \ados{} are also covariant
under translations over half the period: a simultaneous shift of the variables
over $\pi /(2\alpha )$ is equivalent to an interchange of parameters:
\begin{equation}
x_j\ra x_j +\pi /(2\alpha ),\; j=1,\ldots , n \Longleftrightarrow
\mu_0 \leftrightarrow \mu_1 , \;\;\;
\mu_{0}^\prime \leftrightarrow  \mu_{1}^\prime
\end{equation}
(see Eq. \for{VJ} and Eqs. \for{va}, \for{vb}).
\vspace{1ex}

For some purposes it is more convenient to use slightly different expressions
for $\hat{D}_r$. Two such expressions read
\begin{equation}
\hat{D}_{r} =
\sum_{\stackrel{J\subset \{ 1,\ldots , n\} ,\, |J|=r}{\eps_j=\pm 1,\, j\in
J}}\;\;
\sum_{\stackrel{\emptyset\subset J_0\subsetneqq \cdots \subsetneqq J_s=J}
               {0\leq s\leq r} }
(-1)^{s}\prod_{0\leq s^\prime \leq s}
V_{\eps (J_{s^\prime}\setminus J_{s^\prime -1});J_{s^\prime}^c}
 \; e^{-\beta \hat{\theta}_{\eps J_0 } } ,\label{Dr2}
\end{equation}
(with $J_{-1}\equiv\emptyset$) and
\begin{equation}
\hat{D}_{r} =
\sum_{0\leq s\leq r}\;\;
\sum_{\stackrel{J\subset \{ 1,\ldots , n\} ,\, |J|=s}{\eps_j=\pm 1,\, j\in
J}}\;\;
W_{J^c,r-s} V_{\eps J;J^c}
\; e^{-\beta \hat{\theta}_{\eps J}}  , \label{Dr3}
\end{equation}
with
\begin{eqnarray}
W_{I,p}&\equiv&
\sum_{\stackrel{1\leq q \leq p}
               {\eps_i=\pm 1,\, i\in I}}(-1)^{q}
\sum_{
\stackrel{\emptyset\subsetneqq I_1\subsetneqq \cdots \subsetneqq I_{q}\subset
I}
         { |I_q|=p}   } \;
\prod_{1\leq q^{\prime} \leq q}
V_{\eps (I_{q^{\prime}}\setminus I_{q^{\prime} -1});
I\setminus I_{q^{\prime}}}\;\;\;\;\;\;
1\leq p \leq |I|, \nonumber \\
W_{I,0}&\equiv &1 . \label{WJ}
\end{eqnarray}
In Eq. \for{Dr2}, $\hat{D}_r$ is written in terms of the operators
$\exp (-\beta \hat{\theta} )$ rather than the operators
$[\exp (-\beta \hat{\theta} )-1]$. Accordingly, $J_0$ in Eq. \for{Dr2}
is allowed to be empty.
Eq. \for{Dr3} emphasizes the fact that the coefficients of the translator
$\exp (-\beta \hat{\theta}_J )$ in $\hat{D}_r$ consist of
a part $V_{\eps J;J^c}$, which does not commute with the
translator, and a commuting part
$W_{J^c,r-|J|}$.

\vspace{1ex}
\noindent\note From a physical point of view, one can look upon $\hat{D}_r$ as
a
Hamiltonian for an $n$-particle quantum system in dimension one.
The functions of the type $v_a (\pm x_j\pm x_{j^\prime} +2\delta \gamma )$
($\delta =0,1$) and
$v_a(\pm x_j \pm x_k )$ in the coefficients of the \ado{} are in
this interpretation responsible for the interaction between the particles;
the functions $v_b(\pm x_j)$ model an external field.

\subsection{Combinatorial Structure and Parameters}
The increment sets $J_1$, $J_2\setminus J_1,\ldots , J_s\setminus J_{s-1}$
of the increasing sequence \for{part} form the {\em blocks} of
a partition of $J$; the second summation in \for{Dr1} amounts to a sum over all
ordered blocks.
By breaking up $V_{\eps J;K}$ (Eq. \for{VJ}) into three parts
\begin{equation}
V_{\eps J;K}= V^1_{\eps J}\: V^2_{\eps J}\: V^3_{\eps J;K}
\end{equation}
with
\begin{eqnarray}
V^1_{\eps J}&\equiv &\prod_{j\in J} v_b(\eps_j x_j), \label{VJ1} \\
V^2_{\eps J}&\equiv &\prod_{\stackrel{j,j^\prime \in J}{j<j^\prime }}
            v_a(\eps_j x_j+\eps_{j^\prime} x_{j^\prime})\,
		    v_a(\eps_j x_j+\eps_{j^\prime} x_{j^\prime} +2\gamma ),
			\label{VJ2} \\
V^3_{\eps J;K}&\equiv & \prod_{\stackrel{j\in J}{k\in K}}
            v_a(\eps_j x_j +x_k) v_a(\eps_j x_j -x_k),\label{VJ3}
\end{eqnarray}
one can rewrite $\hat{D}_r$ (Eq. \for{Dr1}) as
\begin{eqnarray}
\lefteqn{ \hat{D}_r =
\sum_{\stackrel{J\subset \{ 1,\ldots , n\} ,\, |J|=r}{\eps_j=\pm 1,\, j\in J}}
\left\{ V^1_{\eps J} V^3_{\eps J;J^c}\right. } & &  \label{Dr4} \\
& & \times
\sum_{\stackrel{\emptyset\subsetneqq J_1\subsetneqq\cdots \subsetneqq J_s=J}
               {1\leq s\leq r}} (-1)^{s-1}
\prod_{1\leq s^\prime \leq s}
V^2_{\eps (J_{s^\prime}\setminus J_{s^\prime -1})}
\left.
V^3_{\eps (J_{s^\prime}\setminus J_{s^\prime -1});J\setminus J_{s^\prime}}
\left( e^{-\beta \hat{\theta}_{\eps J_1}} -1 \right) \right\} .\nonumber
\end{eqnarray}
Eqs. \for{Dr1}, \for{Dr2} and \for{Dr3} are more compact than \for{Dr4},
but the latter has the virtue that different parts of the coefficient
can be controlled independently.
The index set $J$ in Eq. \for{Dr4} will be referred to as the {\em cell}.
The first block $J_1$ determines the translator;
this part of the cell will be called the {\em nucleus}.
Notice that: {\em i.} $V^1_{\eps J} V^3_{\eps J;J^c}$ depends on the cell
$J$ but not on its subdivision in blocks;
{\em ii.} the product over $V^2_{\eps (J_{s^\prime}\setminus J_{s^\prime -1})}$
depends on the partition of $J$, but not on the order of the blocks;
{\em iii.} the product over
$V^3_{\eps (J_{s^\prime}\setminus J_{s^\prime -1});J\setminus J_{s^\prime}}$
depends
both on the blocks and on their order.

The parameters $\alpha$ and $\beta$ are scale factors; $\alpha$ determines the
period of
the trigonometric functions and $\beta$ the complex shift of the
translation operators $\exp (\pm \beta \hat{\theta}_j)$.
Both parameters will be taken positive.
The parameters $\mu$, $\mu_\delta$ and $\mu_{\delta^\prime}$ determine the
relative `weight' of $v_a$ and $v_b$ in the coefficients of the \ado{}.
For instance, $v_a\equiv 1$ for $\mu =0$;
therefore, $v_a$ may be omitted
for $\mu = 0$.
The parameters $\mu$, $\mu_\delta$ and $\mu_{\delta^\prime}$ will be assumed to
be non-negative imaginary:
\begin{equation}
\begin{array}{c}
\mu \equiv i\beta g, \;\;\; \mu_\delta \equiv i\beta g_\delta ,\;\;\;
\mu^\prime _\delta \equiv i\beta g^\prime _\delta  ,\;\;\; (\delta =0,1) \\
[1ex]
\alpha ,\beta > 0;\;\;\;\;\;
g,g_\delta ,g^\prime _\delta \geq 0.
\end{array} \label{parres}
\end{equation}
Notice that the above restrictions on the parameters guarantee that:
{\em i.} $\exp (\pm\beta \hat{\theta}_j)$ yields a purely imaginary shift;
{\em ii.} the commuting part of the coefficient, viz.
$W_{J^c,r-s}$ (cf. Eqs. \for{Dr3},\for{WJ}), is real because
$\overline{v_{c}(z)}=v_{c}(-z)$ ($c=a,b$) for $z$ real.

If one picks $\alpha =1/2$, then $\hat{D}_1$ (Eq. \for{D1}) coincides,
up to an irrelevant multiplicative constant,
with Koornwinder's difference operator $D_{\eps_1}$
\cite[Eqs. (5.1)-(5.4)]{koo3}.
The parameters used in \cite{koo3} are related to ours via
\begin{equation}
\label{psubsta}
q=e^{-\beta },\;\;\; t=e^{-\beta g}; \;\;\;\;\;\;
\end{equation}
\begin{equation}
a=e^{-\beta g_0},\;\;\;\;
b=-e^{-\beta g_1},\;\;\;\;
c=e^{-\beta (g^\prime _0 +1/2)},\;\;\;\;
d=-e^{-\beta (g^\prime _1 +1/2)}. \label{psubstb}
\end{equation}

\vspace{1ex}
\noindent\note The parameters $g$, $g_\delta$ and $g^\prime _{\delta}$ can be
physically thought of as being the coupling constants that determine the
strengths of the various interactions. In this interpretation, setting
$g=0$ (i.e. $\mu =0$) yields a system of $n$ particles moving independently
in an external field.

\subsection{$g=0$: Reduction to Rank $n=1$}
By setting $g=0$, the combinatorial structure of $\hat{D}_r$ simplifies
considerably because the coefficients in Eq. \for{Dr4} no longer depend
on the partition of the cell $J$ in blocks
($g=0 \Rightarrow V^2 _{\eps J}=V^3 _{\eps J;K}=1$).
It will be shown next that in this case $\hat{D}_r$ reduces to the $r$th
elementary symmetric function of the following \ados{}:
\begin{equation}
\hat{D}_1 (x_{j})\equiv
\sum_{\eps =\pm 1} v_b (\eps x_j)\left( e^{-\beta\eps\hat{\theta}_j}-1\right) ,
\;\;\;\;\;\; j=1,\ldots , n. \label{D1x}
\end{equation}
(For $n=1$, $\hat{D}_1$ \for{D1} coincides with $\hat{D}_1(x_1)$).

By summation of all terms in $\hat{D}_r$ which correspond to a certain cell $J$
with fixed nucleus $J_1=I$, Eq.~\for{Dr4}
reduces to
\begin{equation}
\hat{D}_r =
\sum_{\stackrel{J\subset \{ 1,\ldots ,n\} ,\, |J|=r}{\eps_j =\pm 1,\, j\in J}}
V^1 _{\eps J}
\sum_{\stackrel{I\subset J}{0\leq |I|\leq r}}
c_{r-|I|}\; e^{-\beta\hat{\theta}_{\eps I}}; \label{reduced}
\end{equation}
where
\begin{equation}
c_0 \equiv 1,\;\;\;\;\;\;\;\;\;\;\;
c_p\equiv \sum_{1\leq s\leq p} (-1)^s N_{p,s}\, ,\;\;\;\;\; p\geq 1,
\label{defc}
\end{equation}
with
\begin{equation}
N_{p,s}= \sum_{\stackrel{p_1+\cdots +p_s=p}{p_j\geq 1}}
\left(
\begin{array}{c}
p \\
p_1 \ldots p_s
\end{array}
\right) ,\;\;\;\;\;\;\; p\geq s\geq 1.
\end{equation}
To verify this,
think of $N_{p,s}$ as the number of ways in which a
collection of $p$ distinct objects can be distributed over $s$ distinct slots
such that every slot is non-empty.

\begin{lem}\hfill

\noindent One has
\begin{equation}
c_p = (-1)^p.\label{cp}
\end{equation}
\end{lem}
\begin{prf}
The above interpretation of $N_{p,s}$ leads to the recursion
relation
\begin{equation}
N_{p,s} = \sum_{s-1\leq q\leq p-1}
\left( \begin{array}{c} p \\ p-q \end{array} \right)
N_{q,s-1}, \;\;\;\;\;\; (p\geq s > 1),\;\;\;\;\;\;\;\;\;\;\;\;\;\;
N_{p,1}= 1. \label{recN}
\end{equation}
Substituting \for{recN} in \for{defc} yields a recursion relation
for $c_p$
\begin{equation}
c_p = -\sum_{0\leq q\leq p-1}
\left( \begin{array}{c} p \\ q \end{array} \right) c_q,
\end{equation}
whose unique solution is \for{cp}.
\end{prf}

\vspace{1ex}
\begin{dfn}\label{elsy}
Let $t_1,\ldots , t_n$ belong to a commutative algebra.
The {\em $r$th elementary symmetric function}
$S_r$ of $t_1,\ldots , t_n$ is defined as
\begin{equation}
S_r(t_1,\ldots , t_n)\equiv
\sum_{\stackrel{ J\subset \{ 1,\ldots , n\} }{|J|=r}} \prod_{j\in J}t_j,
\;\;\;\;\;\; r=1,\ldots , n.
\end{equation}
\end{dfn}
\begin{prp}\label{decoupling}
If $g=0$, then
\begin{equation}
\hat{D}_r =S_r(\hat{D}_1(x_1),\ldots , \hat{D}_1 (x_n)),\;\;\;\;\;\;\;\;
r=1,\ldots , n
\end{equation}
(with $\hat{D}_1 (x_j)$ defined by Eq. \for{D1x}).
\end{prp}
\begin{prf}
Substituting \for{cp} in \for{reduced} yields
\begin{equation}
\hat{D}_r =
\sum_{\stackrel{J\subset \{ 1,\ldots , n\} ,\, |J|=r}{\eps_j =\pm 1,\; j\in J}}
V^1 _{\eps J}
\sum_{\stackrel{I\subset J}{0\leq |I|\leq r}}
(-1)^{r-|I|}e^{-\beta\hat{\theta}_{\eps I}}.\label{red}
\end{equation}
Using Eq.~\for{VJ1}, one completes the proof of the proposition:
\begin{eqnarray}
\hat{D}_r &=&
\sum_{\stackrel{J\subset \{ 1,\ldots , n\} ,\,|J|=r}{\eps_j =\pm 1,\; j\in J}}
\:\prod_{j\in J}  v_b(\eps_j x_j)
\left( e^{-\beta\eps_j \hat{\theta}_j}-1\right)  \\
&=& \sum_{\stackrel{J\subset \{ 1,\ldots , n\}}{|J|=r}}
    \prod_{j\in J}\hat{D}_1 (x_j)\;\;
= S_r(\hat{D}_1(x_1),\ldots , \hat{D}_1(x_n)).
\end{eqnarray}
\end{prf}

\vspace{1ex}
\noindent\note Proposition~\ref{decoupling} is in accordance with the
previously noted fact that for $g=0$ the particles of the
quantum system become independent.

\setcounter{equation}{0}
\addtocontents{toc}{\protect\vspace{-2ex}}
\markright{3 Simultaneous Diagonalization}
\section{Simultaneous Diagonalization}\label{sect3}
In this section it is shown that Koornwinder's polynomials form a basis of
joint eigenfunctions of $\hat{D}_1,\ldots ,\hat{D}_n$. We prove that the
difference operators commute and compute their eigenvalues.
As a result, we obtain an explicit Harish-Chandra-type isomorphism between
the abelian algebra
generated by $\hat{D}_1,\ldots ,\hat{D}_n$ and the symmetric algebra in
$n$ variables.
For convenience, we will put $\alpha =1/2$ from now on.

\subsection{Trigonometric Polynomials}
Let ${\cal A}={\Bbb C}[\exp (ix_1),\ldots , \exp (ix_n)]$ be the algebra of
trigonometric polynomials on the $n$-dimensional torus
\begin{equation}
{\Bbb T}= {\Bbb R}^n / 2\pi {\Bbb Z}^n.\label{torus}
\end{equation}
${\cal A}$ is spanned by the Fourier basis $\{ e^\lambda \}$ with $\lambda$ in
the character lattice ${\cal P}$ of ${\Bbb T}$:
\begin{equation}
e^{\lambda}(x) \equiv  e^{i\sum_{j=1}^{n}\lambda_j x_j},\;\;\;\;\;\;
\lambda \in {\cal P}={\Bbb Z}^n .\label{expo}
\end{equation}
Let $W$ be the (Weyl) group of permutations and sign flips of the variables
$x_1,\ldots  ,x_n$
(so $W\cong S_n \ltimes ({\Bbb Z}_2)^n$).
The subalgebra ${\cal A}^W={\Bbb C}[\cos x_1,\ldots , \cos x_n]^{S_n}$ of
W-invariant
polynomials on ${\Bbb T}$ is spanned by the basis $\{ m_\lambda \}$ of monomial
symmetric
functions
\begin{equation}\label{monomial}
m_{\lambda}(x) = \sum_{\lambda^\prime \in W \lambda} e^{\lambda^\prime}\sim
 \sum_{\lambda^\prime \in S_n \lambda}\;\;
\left( \prod_{1\leq j\leq n} \cos (\lambda^\prime_j x_j) \right) ,\;\;\;\;\;\;
\lambda \in {\cal P^+},
\end{equation}
with $\sim$ denoting proportionality and
${\cal P}^+$ denoting the cone of dominant weights:
\begin{equation}
{\cal P^+} = \{ \lambda \in {\cal P}\; | \;
\lambda_1 \geq \lambda_2 \geq \cdots \geq\lambda_n \geq 0 \; \}
.\label{domcone}
\end{equation}

The lattice ${\cal P}$ can be partially ordered in the
following way:
\begin{dfn}[partial ordering of ${\cal P}$]\label{nat}
\begin{equation}
(\forall \lambda ,\lambda^\prime \in {\cal P}):\;\;\;\;\;\;\;\;\;\;\;
\lambda^\prime \leq \lambda \;\;\;\;\; {\rm iff} \;\;\;\;\;\;
\sum_{1\leq j \leq k} \lambda_j^\prime \leq \sum_{1\leq j\leq k}\lambda_j,
\;\;\;\; {\rm for} \; k=1,\ldots , n.
\end{equation}
\end{dfn}

\vspace{1ex}
The above ordering induces a partial ordering of the monomial basis
$\{ m_\lambda \}_{\lambda \in {\cal P}^+}$.
To each dominant weight $\lambda \in {\cal P}^+$ we associate a finite
dimensional subspace of ${\cal A}^W$ with highest weight $\lambda$:
\begin{equation}
{\cal A}^W_{\lambda} \equiv {\rm span} \{ m_{\lambda^\prime}
\}_{(\lambda^\prime
\in {\cal P}^+,\;\lambda^\prime \leq \lambda )},\;\;\;\;\;\;
\lambda\in {\cal P}^+ .
\end{equation}
Occasionally we will also use the notation
\begin{equation}
|\lambda | \equiv \sum_{1\leq j\leq n} \lambda_j ,\;\;\;\;\;\;\;
 \lambda \in {\cal P}^+.
\end{equation}

\subsection{Triangularity}\label{sectria}
In this subsection it is shown that $\hat{D}_r$ maps the
highest weight spaces ${\cal A}^W_\lambda$ into itself.
\begin{dfn}\label{tria}
A linear operator $\hat{D}:{\cal A}^W \ra {\cal A}^W$ is called
{\em triangular} iff
\begin{equation}
\hat{D} ({\cal A}_\lambda^W) \subset {\cal A}^W_\lambda ,
\;\;\;\;\;\;\;\;\; \forall \lambda \in {\cal P}^+.\label{invhw}
\end{equation}
\end{dfn}
One can rewrite Eq. \for{invhw} in a
more illuminating way:
\begin{equation}
(\forall \lambda \in {\cal P}^+):\;\;\;\;\;\;\;\;\;
\hat{D}\, m_\lambda =
\sum_{\lambda^\prime \in {\cal P}^+,\; \lambda^\prime \leq \lambda}
 [\hat{D}]_{\lambda ,\lambda^\prime}\; m_{\lambda^\prime}\, ,\;\;\;\;\;\;\;\;
{\rm with}\;\; [\hat{D}]_{\lambda ,\lambda^\prime}\in {\Bbb C} \label{triam}
\end{equation}
(i.e. $[ \hat{D} ]_{\lambda ,\lambda^\prime}=0$ if
$\lambda^\prime \nleq \lambda$).
In order to prove that $\hat{D}_r$ is triangular, we first need to verify
that the operator maps ${\cal A}^W$ into itself.
\begin{prp}[invariance of ${\cal A}^W$]\label{inv}
\begin{equation}
\hat{D}_r ({\cal A}^W) \subset {\cal A}^W,\;\;\;\;\;\;\; r=1,\ldots , n.
\end{equation}
\end{prp}
\begin{prf}
Acting with $\hat{D}_r$ Eq. \for{Dr4} on a monomial $m_{\lambda}$
\for{monomial}
yields the following  W-invariant trigonometric function on the torus ${\Bbb
T}$:
\begin{eqnarray}\label{rat}
\lefteqn{(\hat{D}_rm_{\lambda})(x)=
\sum_{\stackrel{J\subset \{ 1,\ldots , n\} ,\, |J|=r}{\eps_j=\pm 1,\, j\in J}}
V^1_{\eps J}\: V^3_{\eps J; J^c}\: \times
\sum_{\stackrel{\emptyset\subsetneqq J_1\subsetneqq\cdots \subsetneqq J_s=J}
               {1\leq s\leq r}  }} & & \\
& &
(-1)^{s-1}\left\{
\prod_{1\leq s^\prime \leq s}
V^2_{\eps (J_{s^\prime}\setminus J_{s^\prime -1})}
\:
V^3_{\eps (J_{s^\prime}\setminus J_{s^\prime -1});
           J\setminus J_{s^\prime}}\;
[ m_\lambda (x+2\gamma e_{\eps J_1})-m_\lambda (x) ] \right\} ,\nonumber
\end{eqnarray}
with
\begin{equation}
e_{\eps J}\equiv \sum_{j\in J} \eps_j e_j
\end{equation}
($\{ e_1,\ldots ,e_n\}$ denotes the standard basis of ${\Bbb R}^n$).
The r.h.s. of \for{rat} is rational in the exponentials $\exp (ix_j)$,
$j=1,\ldots ,n$.
In order to prove the proposition, we need to show that $\hat{D}_r m_\lambda$
\for{rat} is actually a polynomial in $\exp (ix_1),\ldots ,\exp (ix_n)$.
Since the r.h.s. of \for{rat} is symmetric in $x_1,\ldots , x_n$ it suffices to
verify that $\hat{D}_r m_\lambda$, viewed as a function of $x_1$, is free of
poles.

As a function of $x_1$, the terms in \for{rat} may have poles
caused by zeros in the denominators of the coefficients of
the \ado{}. These poles are located at
(cf. Eqs. \for{VJ1}-\for{VJ3} and \for{va}, \for{vb}):
\begin{equation}
\begin{array}{ccll}
x_1 &=& 0  &     {\rm mod}\; \pi , \\
    &=& \pm \gamma &   {\rm mod}\; \pi , \\
	&=& \pm x_j &   {\rm mod}\; 2\pi ,\;\; j=2,\ldots , n,  \\
    &=& \pm x_j\pm 2\gamma & {\rm mod}\; 2\pi ,\;\; j=2,\ldots , n.
\end{array}\label{locpol}
\end{equation}
{}From now on the parameters
$\gamma$, $\mu$, $\mu_\delta$, $\mu^\prime_\delta$
($\delta =0,1$) and
the remaining variables $x_2,\ldots , x_n$ are fixed in general position.
Specifically, we choose these parameters and variables such that
the poles in the terms of \for{rat} are simple.

The residue at $x_1=0$ vanishes because \for{rat}
is even in $x_1$;
the residue at $x_1=x_j$ vanishes because
\for{rat} is invariant under an interchange of the variables $x_1$ and $x_j$.
Furthermore, because $\hat{D}m_\lambda (x)$ is even in $x_j$,
$j=1,\ldots , n$ and covariant
under translations over half the period
(cf. Remark {\em vi}, Section~\ref{ados}):
\begin{equation}
x_j \ra x_j +\pi ,\; j=1,\ldots ,n \Longleftrightarrow
\left\{ \begin{array}{l}
       \mu_0 \leftrightarrow \mu_1,\;\;\;
	   \mu^\prime_0 \leftrightarrow \mu^\prime_1 \\
	   m_{\lambda}(x) \ra (-1)^{|\lambda |} m_\lambda (x), \label{cov}
	   \end{array}  \right.
\end{equation}
(with $|\lambda |= \sum_{j=1}^n \lambda_j$),
we need only show that the total residue of \for{rat}
vanishes at e.g.:
\begin{equation}
\begin{array}{cclll}
x_1 &=& - \gamma &    & {\rm (type\; I)}, \\
 \;	&=& - x_j-2\gamma  &   j=2,\ldots , n &
	                                {\rm (type\; II)}.
\end{array}
\end{equation}

type I: $x_1=-\gamma$.\\
The only terms in the r.h.s. of \for{rat} that contribute to the residue at
$x_1=-\gamma$ are those corresponding to cells $J$ with $1\in J$ and
$\eps_1 =+1$ (recall \for{vb} and \for{VJ1} to check this).
Fix a cell $J$ and choose a configuration of signs
$\eps_j\in \{+1,-1\} $, $j\in J$. It suffices to show that the total residue at
$x_1=-\gamma$ in the sum of all terms of \for{rat} corresponding to
this fixed cell $J$, with the signs prescribed, is zero.
One may assume that all signs $\eps_j$,
$j\in J$ are positive; this is because the general situation can be
reduced to the case with $\eps_j=1$ by appropriate flipping
of signs of the variables $x_j$, $j\in J$.
We prove in Appendix~A (Lemma~\ref{pole1}) that the sum of terms
in \for{rat} corresponding to a fixed cell $J$ with all signs positive
is indeed regular at $x_1=-\gamma$.

type II: $x_1=-x_j-2\gamma$.\\
The proof of this case is very similar to the previous one.
The only terms in \for{rat} that contribute to the residue at
$x_1=-x_j-2\gamma$ are those with $1,j\in J$ and $\eps_1=\eps_j=+1$
(cf. \for{va} and \for{VJ2}).
Lemma~\ref{pole2} of Appendix~A states that the sum of all terms in
\for{rat} that correspond to a fixed cell $J$
with the signs $\eps_j$, $j\in J$, being $+1$, is regular at
$x_1=-x_j-2\gamma$.
Again the general case (corresponding to an arbitrary configuration
of signs $\eps_j=\pm 1$) can be obtained by an appropriate flipping of the
signs
of $x_j$, $j\in J$.

We conclude that $\hat{D}_r m_\lambda$ is a $W$-invariant rational
function on the torus ${\Bbb T}$ without poles.
Consequently, $\hat{D}_r m_\lambda$ must be a polynomial in
${\cal A}^W$, which completes the proof of the proposition.
\end{prf}

Proposition~\ref{inv} says that $\hat{D}_r m_\lambda$ is a W-invariant
trigonometric polynomial on ${\Bbb T}$, i.e. it must be a finite linear
combination of
monomial symmetric functions:
\begin{equation}
(\forall \lambda \in {\cal P}^+):\;\;\;\;\;\;\;\;\;\;
\hat{D_r}m_\lambda =
\sum_{\lambda^\prime \in {\cal P}^+_{\lambda ,r} }
 [\hat{D}_r]_{\lambda ,\lambda^\prime}m_{\lambda^\prime},
\;\;\;\;\;\;\;\; |{\cal P}^+_{\lambda ,r} |<\infty , \label{expand}
\end{equation}
with
\begin{equation}
{\cal P}^+_{\lambda ,r} \equiv \{ \lambda^\prime \in {\cal P}^+ \; |\;
[\hat{D}_r]_{\lambda ,\lambda^\prime}\neq 0 \} .\label{nonzerot}
\end{equation}
In order for $\hat{D}_r$ to be triangular one must have
\begin{equation}
\lambda^\prime \in {\cal P}^+_{\lambda ,r}\Longrightarrow
\lambda^\prime \leq \lambda \label{trian}
\end{equation}
(cf. Eq. \for{triam}).
We shall prove this property by studying the asymptotics of
$(\hat{D}_r\, m_{\lambda})(x)$ for Im~$x_j\ra -\infty$.
The following limits will be useful (cf. Eqs. \for{va}, \for{vb} and
\for{parres}):
\begin{eqnarray}
\lim_{R \ra \infty} v_a (z+i\eps R ) &=&e^{\eps\beta g/2},\label{asva}\\
\lim_{ R \ra \infty} v_b (z+i\eps R ) &=&
e^{\eps\beta (g_0+g_1+g_0^\prime +g_1^\prime )/2} \label{asvb}
\end{eqnarray}
(with $\eps =\pm 1$).
\begin{prp}[triangularity]\label{Dtri}
\begin{equation}
(\forall \lambda \in {\cal P}^+):\;\;\;\;\;
\hat{D}_r ({\cal A}^W_\lambda ) \subset {\cal A}^W_\lambda ,\;\;\;\;\;\;
r=1,\ldots , n.
\end{equation}
\end{prp}
\begin{prf}
Fix an $r\in \{ 1,\ldots , n\} $ and $\lambda \in {\cal P}^+$.
Let $\omega_k\equiv \sum_{1\leq j\leq k} e_j$ (the $k$th fundamental weight)
and introduce (cf.  \for{nonzerot})
\begin{equation}
M_{\lambda ,r;k}\equiv
{\rm max} \{ (\lambda^\prime ,\omega_k) \;\; |\;
 \lambda^\prime \in {\cal P}^+_{\lambda ,r}\;   \} .
\end{equation}
To derive a contradiction, assume
\for{trian} does not hold;
i.e. assume that there exists a $k\in \{ 1,\ldots , n\} $ such that
\begin{equation}
 M_{\lambda , r ;k}> (\lambda ,\omega_k)
\left( =\sum_{1\leq j\leq k} \lambda_j \right) .\label{ineq}
\end{equation}
Now it is easy to verify the asymptotics
\begin{equation}
m_{\lambda^\prime}( x-i R\omega_k) \sim
e^{ R\, (\lambda^\prime ,\omega_k)}
\sum_{\lambda^{\prime\prime} \in W_{\lambda^\prime ;k} (\lambda^\prime )}\;
e^{\lambda^{\prime \prime}}(x),\;\;\;\;\;\;\;\;  R \ra +\infty ,\label{monas}
\end{equation}
with
\begin{equation}
W_{\lambda^\prime ;k} \equiv
\{ w\in W \; |\; (w\lambda^\prime ,\omega_k)=(\lambda^\prime ,\omega_k)
\;\;   \} ,
\end{equation}
so using \for{expand} we obtain
\begin{equation}
\lim_{ R \ra \infty}
e^{- R M_{\lambda ,r;k}}(\hat{D}_r
      m_{\lambda})( x-i R \omega_k) =
\sum_{ \stackrel{\lambda^\prime \in {\cal P}^+_{\lambda ,r}}
                {(\lambda^\prime ,\omega_k) = M_{\lambda ,r;k}} }
[\hat{D}_r]_{\lambda ,\lambda^\prime}\:
\left( \sum_{\lambda^{\prime\prime} \in W_{\lambda^\prime ;k}(\lambda^\prime
)}\:
 e^{\lambda^{\prime\prime}}(x)
\right) .\label{lim1}
\end{equation}

On the other hand, Eq.~\for{monas} combined with the limits \for{asva} and
\for{asvb} entails the following asymptotics for \for{rat}:
\begin{equation}
(\hat{D}_rm_\lambda )(x-i R\omega_k) =
O( e^{ R\, (\lambda ,\omega_k) } ),\;\;\;\;\;\;\;
 R \ra \infty .
\end{equation}
Consequently,
\begin{equation}
\lim_{ R \ra \infty}
e^{- R M_{\lambda ,r;k}}(\hat{D}_r
m_{\lambda})(x-i R\omega_k)=0\label{lim2}
\end{equation}
(because of inequality \for{ineq}).

The matrix elements $[\hat{D}_r]_{\lambda ,\lambda^\prime}$ in
\for{lim1} are non-zero (by definition),
and the exponentials $e^{\lambda^{\prime\prime}}(x)$ \for{expo} corresponding
to
different weights $\lambda^{\prime\prime } \in {\cal P}$ are linearly
independent.
Hence, by comparing the r.h.s. of
Eqs. \for{lim1} and \for{lim2} one arrives at the desired contradiction.
\end{prf}

\subsection{The Spectrum}
By extending $\leq$ (Definition~\ref{nat})
to a linear ordering of ${\cal P}^+$ (take e.g. the Lexicographical
ordering), it is easy to see that the triangularity of the \ado{}:
\begin{equation}
\hat{D}_r m_\lambda =
\sum_{\lambda^\prime \in {\cal P}^+,\, \lambda^\prime \leq \lambda}
[\hat{D}_r]_{\lambda ,\lambda^\prime}\; m_{\lambda^\prime}\, ,
\;\;\;\;\;\;\;\;\;\;\; \lambda \in {\cal P}^+, \label{expand2}
\end{equation}
has as consequence that the elements on the diagonal of the matrix
$[\hat{D}_r]_{\lambda ,\lambda^\prime}$ are the eigenvalues of the operator
$\hat{D}_r: {\cal A}^W \ra {\cal A}^W$.
(It will become clear in Section~\ref{diacom} that in our case
these eigenvalues are semisimple).
The purpose of the present subsection is to compute
$[\hat{D}_r]_{\lambda ,\lambda}$.

Let $y\in {\Bbb R}^n$ be a fixed vector subject to the condition
\begin{equation}
y_1> y_2 >\cdots > y_n >0. \label{chamber}
\end{equation}
By combining the asymptotics (cf. \for{monas})
\begin{equation}
m_{\lambda^\prime} (i R y) \sim e^{ R \sum_{j=1}^n \lambda^\prime_j y_j}
\;\;\;\;\;\;\;\;  R \ra \infty ,\label{aschamber}
\end{equation}
with Eq. \for{expand2}, one finds
(use $\lambda^\prime < \lambda \Rightarrow
\sum_{j=1}^n \lambda_j^\prime y_j < \sum_{j=1}^n \lambda_j y_j$)
\begin{equation}
[\hat{D}_r]_{\lambda ,\lambda}=
\lim_{ R \ra \infty} e^{- R \sum_{j=1}^n \lambda_j y_j}
(\hat{D}_r m_\lambda )(i R y). \label{limdiag}
\end{equation}
In the next proposition, we will evaluate this limit.
\begin{prp}[eigenvalues]\label{Deigv}\hfill

\noindent One has
\begin{equation}
[\hat{D}_r]_{\lambda ,\lambda}=2^r\:
E_{r,n} ( \ch\beta (\lambda_1+\rho_1),\ldots , \ch\beta (\lambda_n+\rho_n);
          \ch\beta \rho_r,\ldots , \ch\beta \rho_n ),\label{ev}
\end{equation}
with
\begin{eqnarray}
\lefteqn{
E_{r,n} (t_1,\ldots , t_n;p_r,\ldots , p_n)\equiv } & & \nonumber \\
& & \sum_{0\leq s\leq r} (-1)^{r+s}
\left( \sum_{\stackrel{J\subset \{ 1,\ldots , n\} }{|J|=s}}\;
       \prod _{j\in J} t_j \right)
\left( \sum_{r\leq i_1\leq\cdots \leq i_{r-s}\leq n} \;
       p_{i_1}\cdots p_{i_{r-s}} \right)
\end{eqnarray}
and
\begin{equation}
\rho_j \equiv (n-j)g + (g_0 +g_1 +g_0^\prime +g_1^\prime )/2,
\;\;\;\;\;\;\;\; j=1,\ldots , n.\label{roo}
\end{equation}
\end{prp}
\begin{prf}
Using Eq.~\for{Dr3} we obtain
\begin{equation}
(\hat{D}_{r}m_\lambda )(x) =
\sum_{0\leq s\leq r}\;\;
\sum_{\stackrel{J\subset \{ 1,\ldots , n\} ,\, |J|=s}{\eps_j=\pm 1,\, j\in
J}}\;\;
W_{J^c,r-s} V_{\eps J;J^c}\; m_\lambda (x+i\beta e_{\eps J}),\label{Dract}
\end{equation}
with $V_{\eps J;K}$ and $W_{I,p}$ defined by Eqs. \for{VJ} and \for{WJ},
respectively.
In order to compute the limit \for{limdiag}, we first derive some preliminary
asymptotics:

{\em i.} $m_\lambda$ (cf. Eq. \for{monas}):
\begin{equation}
m_\lambda (x+i\beta e_{\eps J})|_{x=i R y} \sim
e^{ R \sum_{j=1}^n \lambda_j y_j }\;
e^{\beta \sum_{j\in J} \eps_j \lambda_j },\;\;\;\;\;\;\;  R \ra \infty .
\label{as1}
\end{equation}

{\em ii.} $V_{\eps J;J^c}\:$:\\
{}From \for{asva}, \for{asvb} one deduces
\begin{equation}
\lim_{ R \ra \infty} V_{\eps J;J^c}|_{x=i R y} =
e^{\beta [ gN_{\eps J}+ (g_0+g_1+g_0^\prime +g_1^\prime )M_{\eps J}]/2},
\end{equation}
with
\begin{eqnarray}
M_{\eps J}&=&\left| \{ j\in J \; |\; \eps_j =+1 \} \right| -
             \left| \{ j\in J \; |\; \eps_j =-1 \} \right| \nonumber \\
	  &=&\sum_{j\in J} \eps_j
\end{eqnarray}
and
\begin{eqnarray}
N_{\eps J} &=& 2 \left| \{ j,k \in J \; |\;  j<k,\; \eps_j=+1 \} \right| \;\,
\;\;\;\;\;\;\;  -2 \left| \{ j,k \in J \; |\;  j<k,\; \eps_j=-1 \} \right| +
\nonumber \\
   & & 2 \left| \{ j\in J,\; k \in J^c \; |\;  j<k,\; \eps_j=+1 \} \right| -
  2\left| \{ j\in J,\; k \in J^c \; |\;  j<k,\; \eps_j=-1 \} \right|
\nonumber\\
   &=& 2\left| \{ j\in J,\; k=1,\ldots , n\; |\;  j<k,\; \eps_j=+1 \} \right|
-\nonumber \\
   & &  2 \left| \{ j\in J,\; k=1,\ldots , n\; |\;  j<k,\; \eps_j=-1 \} \right|
\nonumber \\
   &=& 2\sum_{j\in J} (n-j)\eps_j .
\end{eqnarray}
Consequently,
\begin{equation}
\lim_{ R \ra \infty} V_{\eps J;J^c}|_{x=i R y} =
e^{\beta \sum_{j\in J} \eps_j \rho_j}\label{as2}
\end{equation}
with $\rho_j$ defined by \for{roo}.

{\em iii.} $W_{I,p}\:$:\\
Using \for{WJ} and \for{as2} it is not hard to see that
$\lim_{ R \ra \infty} W_{I,p}|_{x=i R y}$ exists and depends only
on the cardinality of $I$ and on $p$ (but not on the number of variables $n$).
We define
\begin{equation}
F_{|I|,p} \equiv 2^{-p}
\lim_{ R \ra \infty} W_{I,p}|_{x=i R y}.\label{as3}
\end{equation}
Notice that (cf. Eq. \for{WJ})
\begin{equation}
F_{m,0}=1,\;\;\;\;\;\;\;\;\; m =0,\ldots , n.\label{init}
\end{equation}

After these preliminaries, we are now ready to compute limit \for{limdiag}.
Substituting \for{Dract} in \for{limdiag}, and making use of \for{as1},
\for{as2}
and \for{as3}, we obtain
\begin{equation}
[\hat{D}_r]_{\lambda ,\lambda}=2^r
\sum_{\stackrel{J\subset \{ 1,\ldots , n\} ,\, |J|=s}{ 0\leq s\leq r}}
\left( \prod _{j\in J} \ch \beta (\lambda_j +\rho_{j}) \right)
F_{n-s,r-s}.\label{evF}
\end{equation}
It remains to calculate $F_{m,p}$, $1\leq p\leq m\leq n$.
Acting with $\hat{D}_r$ on a constant function yields zero
(see Eq. \for{Dr1}); one obtains, therefore, the following relations
for $F_{m,p}$ after setting $\lambda =0$:
\begin{equation}
\sum_{\stackrel{J\subset \{ 1,\ldots , n\} ,\, |J|=s}{ 0\leq s\leq r}}
\left( \prod _{j\in J} \ch (\beta \rho_j) \right) F_{n-s,r-s}=0,
\;\;\;\;\;\;\; 1\leq r\leq  n.\label{recF}
\end{equation}
For a fixed number of variables $n$, this yields $n$ equations.
However,
for $n^\prime < n$ the same coefficients $F_{m,p}$ occur,
now with $0\leq p\leq m\leq n^\prime$.
Collecting the relations \for{recF} for $n^\prime =1,\ldots ,n$, and making
use of Eq.~\for{init} results in a
linear system of $n(n+1)/2$ equations in the $n(n+1)/2$ variables $F_{m,p}$,
$1\leq p \leq m \leq n$.

In order to solve this system it is convenient to make the substitution
$\rho_{j}=\tilde{\rho}_{n+1-j}$. Notice that $\tilde{\rho}_j$
(unlike $\rho_j$) does not depend on the number of variables $n$.
After this substitution the $n(n+1)/2$ equations become:
\begin{equation}
\sum_{\stackrel{J\subset \{ 1,\ldots , n^\prime \} ,\, |J|=s}{ 0\leq s\leq r}}
\left( \prod _{j\in J} \ch (\beta \tilde{\rho}_j) \right) F_{n^\prime
-s,r-s}=0,
\;\;\;\;\;\;\; 1\leq r\leq n^\prime \leq n,\label{sst}
\end{equation}
with condition \for{init}.
In Lemma~\ref{reclem1} of Appendix~B it is shown that this linear system
has a unique solution:
\begin{equation}
F_{m,p} = (-1)^{p}
\sum_{1\leq i_1\leq i_2\leq \cdots \leq i_{p}\leq m+1-p}
\; \ch (\beta\tilde{\rho}_{i_1})\cdots \ch (\beta \tilde{\rho}_{i_{p}})
\label{sol}
\end{equation}
Substituting \for{sol} in \for{evF} and using $\tilde{\rho}_j=\rho_{n+1-j}$
now yields the expressions \for{ev}-\for{roo}.
\end{prf}

For $g=0$ all components of the vector $\rho$ (Eq. \for{roo}) are
equal:
\begin{equation}
\rho_1,\rho_2,\ldots ,\rho_n=
(g_0+g_1+g_0^\prime +g_1^\prime)/2 \equiv \overline{\rho}.
\end{equation}
Then, one can rewrite the above expressions for the eigenvalues in
terms of elementary symmetric functions (see the remark following
Lemma~\ref{reclem2} of Appendix~B. to check this):
\begin{eqnarray}
[\hat{D}_r]_{\lambda ,\lambda}&=&
2^r \sum_{\stackrel{J\subset \{ 1,\ldots , n\}}{|J|=r}}
\prod_{j\in J} \left( \ch\beta (\lambda_j +\overline{\rho})-
                      \ch\beta \overline{\rho} \right) \nonumber \\
&=& 2^r\: S_r\left(
\ch\beta (\lambda_1 +\overline{\rho})-\ch\beta \overline{\rho},\ldots ,
\ch\beta (\lambda_n +\overline{\rho})-\ch\beta \overline{\rho} \right) .
\label{specdec}
\end{eqnarray}
This equation is in agreement with Proposition~\ref{decoupling}.

\subsection{Symmetry}
In Ref.~\cite{koo3} the following weight function on
the torus ${\Bbb T}$ was introduced:
\begin{equation}
\Delta (x) \equiv \prod_{1\leq j< k\leq n} d_a(x_j+x_k)\, d_a(x_j-x_k)
                  \prod_{1\leq j\leq n} d_b(x_j) \label{densityb}
\end{equation}
with
\begin{eqnarray}
d_{c}(z)&=& d_{c}^+(z)d_{c}^+(-z),\;\;\;\;\;\;\;\;\;\;\;\;\;\;\; c=a,b
\label{dab} \\
 d_{a}^+(z)&\equiv& \frac{ (e^{iz};e^{-\beta})_\infty }
                        { (e^{-\beta g}e^{iz};e^{-\beta})_\infty },\label{da+}
\;\;\;\;\;\;\;\;\;\; \\
d_{b}^+(z)&\equiv& \frac{ (e^{2iz};e^{-\beta})_\infty }
   { (e^{-\beta g_0}e^{iz},-e^{-\beta g_1}e^{iz},
       e^{-\beta (g_0^\prime +1/2)}e^{iz},
	  -e^{-\beta (g_1^\prime +1/2)}e^{iz}; e^{-\beta})_\infty }.\label{db+}
\end{eqnarray}
The so-called q-shifted factorials are defined in the usual way:
\begin{equation}
(a;q)_\infty \equiv \prod_{m=1}^\infty (1-a q^m),\;\;\;\;\;\;\;
(a_1,\ldots , a_r;q)_\infty \equiv \prod_{s=1}^r (a_s;q)_\infty .
\end{equation}
Notice that the conditions on our parameters, viz. \for{parres}, guarantee that
the infinite products in Eqs. \for{da+} and \for{db+} converge. Recall that in
order to compare our formulas with those of \cite{koo3}, one has to
reparametrize
according to Eqs. \for{psubsta} and \for{psubstb}.

Let $L^2_W({\Bbb T},\Delta dx)$ be the space of W-invariant functions on
${\Bbb T}$ that
are square integrable with respect to the measure $\Delta dx$.
We define
\begin{equation}
\langle f ,g \rangle_\Delta \equiv
\int_{\Bbb T} f\overline{g} \Delta dx,\;\;\;\;\;\;\;\; f,g\in L^2_W({\Bbb
T},\Delta dx).
\label{ip}
\end{equation}
The space of W-invariant polynomials  ${\cal A}^W$ is a dense subspace of
$L^2_W({\Bbb T},\Delta dx)$.
The purpose of the present section is to show that the \ados{}
$\hat{D}_1,\ldots ,\hat{D}_n$
are symmetric with respect to
$\langle \cdot ,\cdot \rangle_\Delta$.
We need the following lemma:

\begin{lem}
Let $\alpha =1/2$ and $z\in {\Bbb R}$; furthermore, let the parameters be
subject to  condition
\for{parres}. Then the functions $d_{c}^{(+)}$ ($c=a,b$) satisfy the following
first order difference equations:
\renewcommand{\theenumi}{\roman{enumi}}
\begin{enumerate}
\item $d_{c}^{+}$:
\begin{eqnarray}
d_{a}^+(z+i\beta )&=&e^{-\beta g/2}\: v_a(z)\: d_a^+(z), \label{d1}\\
d_{b}^+(z+i\beta )&=&e^{-\beta (g_0+g_1+g_0^\prime +g_1^\prime )/2}\:
v_b(z)\: d_b^+(z); \label{d2}
\end{eqnarray}
\item $d_{c}$:
\begin{equation}
v_{c}(z-i\beta )\: d_{c}(z-i\beta )=\overline{v_{c}(z)}\: d_{c}(z),
\;\;\;\;\;\;\; c=a,b.\label{d3}
\end{equation}
\end{enumerate}
\end{lem}
\begin{prf}
{\em i.} Eq. \for{d1} is an immediate consequence of definition \for{da+}:
\begin{equation}
d_{a}^+(z+i\beta )/d_a^+(z)=\frac{(1-e^{-\beta g}e^{iz})}{(1-e^{iz})}=
e^{-\beta g/2}\: v_a(z).
\end{equation}
Eq. \for{d2} can be reduced to the former case by observing that $d_b^+(z)$
factorizes:
\begin{eqnarray}
d_b^+(z)&=&\frac{ (e^{iz};e^{-\beta})_\infty }
              { (e^{-\beta g_0}e^{iz};e^{-\beta})_\infty }
         \frac{ (-e^{iz};e^{-\beta})_\infty }
              { (-e^{-\beta g_1}e^{iz};e^{-\beta})_\infty } \nonumber \\
& & \times  \frac{ (e^{-\beta /2}e^{iz};e^{-\beta})_\infty }
              { (e^{-\beta (g_0^\prime +1/2)}e^{iz};e^{-\beta})_\infty }
         \frac{ (-e^{-\beta /2}e^{iz};e^{-\beta})_\infty }
       { (-e^{-\beta (g_1^\prime +1/2)}e^{iz};e^{-\beta})_\infty }.\label{db+f}
\end{eqnarray}.

\noindent {\em ii.} Using Eqs. \for{d1} or \for{d2}, respectively
(in step * below), one derives Eq. \for{d3}:
\begin{eqnarray}
v_{c}(z-i\beta )\: d_{c}(z-i\beta )&=&
v_{c}(z-i\beta )\: d^+_{c}(z-i\beta )\cdot d^+_{c}(-z+i\beta ) \nonumber \\
&\stackrel{*}{=}&
d^+_{c}(z)\cdot v_{c}(-z)\: d^+_{c}(-z) \nonumber \\
&=& \overline{v_{c}(z)}\: d_{c}(z).
\end{eqnarray}
\end{prf}
Part {\em ii.} of the above lemma leads to the following useful difference
equation:
\begin{cor}\label{fsym}One has
\begin{equation}
e^{\beta \hat{\theta}_{\eps J}}(V_{\eps J;J^c}\Delta )=
\overline{V_{\eps J;J^c}}\Delta .\label{fsymeq}
\end{equation}
\end{cor}

We now arrive at the main result of this subsection namely the symmetry
of $\hat{D}_r$, whose proof hinges on relation \for{fsymeq}.
\begin{prp}[symmetry]\label{Dsym}
\begin{equation}
\langle \hat{D}_r m_\lambda ,m_{\lambda^\prime} \rangle_\Delta
=\langle m_\lambda ,\hat{D}_r m_{\lambda^\prime} \rangle_\Delta ,
\;\;\;\;\;\;\;\;\;\forall\: \lambda ,\lambda^\prime \in {\cal P}^+.\label{sym}
\end{equation}
\end{prp}
\begin{prf}
First consider the following contour integral
\begin{equation}
\oint_{C_j}W_{J^c,r-s}\: V_{\eps J;J^c}\:
( e^{-\beta\hat{\theta}_{\eps J}}m_\lambda )\; \overline{m_{\lambda^\prime}}
\; \Delta\, dx_j\label{int1}
\end{equation}
with $V_{\eps J;K}$ and $W_{I,p}$ as in \for{VJ} and \for{WJ}, respectively
and $j\in J$.
The integration takes place over the closed contour
\begin{equation}
C_j = [-\pi ,\pi ]\cup
      [\pi ,\pi -i\eps_j\beta ]\cup
      [\pi -i\eps_j\beta ,-\pi -i\eps_j\beta ]\cup
	  [-\pi -i\eps_j\beta ,-\pi ].
\end{equation}
Let all parameters and the variables $x_k$, $k\neq j$, be fixed in general
position. A priori the integrand has simple poles inside $C_j$
due to zeros in the denominators of $V_{\eps J;J^c}$ and $\Delta (x)$.
However, one easily verifies that any of these poles in $V_{\eps J;J^c}$ is
compensated by a zero in $\Delta (x)$; similarly, poles inside $C_j$
due to $\Delta (x)$ are compensated by a zero in $V_{\eps J;J^c}$.
Consequently, integral \for{int1} vanishes because of Cauchy's theorem.
Furthermore, the contributions to \for{int1} which are due to the paths
$[\pi ,\pi -i\eps_j\beta ]$ and $[-\pi -i\eps_j\beta ,-\pi ]$ respectively,
cancel each other because the integrand is periodic in $x_j$ with period
$2\pi$.
The upshot is that, when integrating the integrand of \for{int1} over
$[-\pi ,\pi ]$, one may deform the integration path to
$[-\pi -i\eps_j\beta ,\pi -i\eps_j\beta ]$ without changing the value of the
integral.

Armed with this conclusion and Eq. \for{fsymeq}, we are now ready to
prove \for{sym}.
Its l.h.s. can be written
\begin{equation}
\langle \hat{D}_r m_\lambda ,m_{\lambda^\prime} \rangle_\Delta
\stackrel{{\rm Eq.} \for{Dr3}}{=}
\sum_{\stackrel{|J|=s,\, 0\leq s\leq r}{\eps_j=\pm 1,\, j\in J}}
\int_{\Bbb T} W_{J^c,r-s}\: V_{\eps J;J^c}\:
( e^{-\beta\hat{\theta}_{\eps J}}m_\lambda )\; \overline{m_{\lambda^\prime}}
\; \Delta \, dx.\label{int2}
\end{equation}
Deformation of the integration paths of $x_j$, $j\in J$, from $[-\pi ,\pi]$ to
$[-\pi -i\eps_j\beta ,\pi -i\eps_j \beta ]$, followed by a change
of variables $x_j \ra x_j -i\beta \eps_j$, $j\in J$, yields
\begin{equation}
\langle \hat{D}_r m_\lambda ,m_{\lambda^\prime} \rangle_\Delta
= \sum_{\stackrel{|J|=s,\, 0\leq s\leq r}{\eps_j=\pm 1,\, j\in J}}
\int_{\Bbb T} W_{J^c,r-s}\: m_\lambda \;
( e^{\beta\hat{\theta}_{\eps J}}\; \overline{m_{\lambda^\prime}}\:
V_{\eps J;J^c}\: \Delta )dx.\label{int3}
\end{equation}
Using Corollary~\ref{fsym} and the fact that $W_{J^c,r-s}$ is real
(for parameters subject to \for{parres}) entails
\begin{eqnarray}
\langle \hat{D}_r m_\lambda ,m_{\lambda^\prime} \rangle_\Delta &=&
\sum_{\stackrel{|J|=s,\, 0\leq s\leq r}{\eps_j=\pm 1,\, j\in J}}\int_{\Bbb T}
m_\lambda\; \overline{ W_{J^c,r-s}\: V_{\eps J;J^c} \;
(e^{-\beta \hat{\theta}_{\eps J}}\: m_{\lambda^\prime}) }\;\Delta \, dx
\nonumber \\
&=& \langle m_\lambda ,\hat{D}_r m_{\lambda^\prime} \rangle_\Delta .
\end{eqnarray}
\end{prf}

\subsection{Diagonalization and Commutativity}\label{diacom}
If $g=g_0=g_1=g_0^\prime =g_1^\prime =0$, then $d_a =d_b =1$
(see Eqs. \for{dab},\for{da+} and \for{db+f}), and thus $\langle\cdot
,\cdot\rangle_\Delta$
reduces to the inner product on ${\Bbb T}$ with respect to Lebesgue
measure ($\Delta =1$). The basis of monomials
$\{ m_{\lambda}\}_{\lambda\in {\cal P}^+}$ is an orthogonal basis of
$L^2_W({\Bbb T},\, dx)$. For arbitrary parameters however, the orthogonality of
the monomials with respect to $\langle\cdot ,\cdot\rangle_\Delta$ no longer
holds.
By subtracting from $m_\lambda$ the orthogonal projection of $m_\lambda$ onto
span$\{ m_{\lambda^\prime}\}_{
(\lambda^\prime\in {\cal P}^+,\,\lambda^\prime <\lambda )}
\subset {\cal A}^W_\lambda$, one obtains an
alternative basis $\{ p_\lambda \}_{\lambda \in {\cal P}^+}$ of ${\cal A}^W$.
This is the basis of Koornwinder polynomials.
\begin{dfn}[Koornwinder polynomials]\label{koopol}\hfill

\noindent {\em Koornwinder's polynomial} $p_\lambda \in {\cal A}^W_\lambda$
is defined by the conditions
\begin{equation}
p_\lambda = m_\lambda +
\sum_{\lambda^\prime \in {\cal P}^+,\, \lambda^\prime < \lambda}
c_{\lambda ,\lambda^\prime}\: m_{\lambda^\prime}  \label{kp1}
\end{equation}
and
\begin{equation}
\langle p_\lambda ,m_{\lambda^\prime}\rangle_\Delta =0,
\;\;\;\;\;\; \forall \lambda^\prime \in {\cal P}^+,\;
\lambda^\prime <\lambda .\label{kp2}
\end{equation}
\end{dfn}
We now prove that $\{ p_\lambda \}_{\lambda \in {\cal P}^+}$
is a basis of joint eigenfunctions of $\hat{D}_1,\ldots ,\hat{D}_n$.
For convenience, the notation for the eigenvalues Eq. \for{ev} is sometimes
abbreviated by putting
\begin{equation}
E_{r,n}(\ch\beta\theta_1,\ldots ,\ch\beta\theta_n;
        \ch\beta\rho_r,\ldots ,\ch\beta\rho_n)\longrightarrow
E_{r,n}(\theta),\;\;\;\;\; (\theta\in {\Bbb R}^n). \label{ern}
\end{equation}

\begin{thm}[eigenfunctions]\label{efunc}
\begin{equation}
\hat{D}_r\: p_\lambda = E_{r,n}(\lambda +\rho )\: p_\lambda ,
\;\;\;\;\;\;\;\;\;\;\;\; \forall \lambda \in {\cal P}^+\label{Dreveq}
\end{equation}
(with $r=1,\ldots ,n$ and $\rho =(\rho_1,\ldots , \rho_n)$, see Eq.~\for{roo}).
\end{thm}
\begin{prf}
It follows from \for{kp1} and Proposition~\ref{Dtri} that
\begin{equation}
(\hat{D}_r -[\hat{D}_r]_{\lambda ,\lambda})\: p_{\lambda} \in \,
{\cal A}^W_{\lambda}.\label{inco}
\end{equation}
On the other hand, one has (use the propositions \ref{Dsym}, \ref{Dtri} and
Eq.~\for{kp2})
\begin{equation}
\langle (\hat{D}_r -[\hat{D}_r]_{\lambda ,\lambda})\: p_{\lambda},
m_{\lambda^\prime} \rangle_\Delta
\stackrel{{\rm symmetry}}{=}
\langle p_\lambda ,
(\hat{D}_r -[\hat{D}_r]_{\lambda ,\lambda})\: m_{\lambda^\prime} \rangle_\Delta
\stackrel{\for{expand2},\for{kp2}}{=} 0,
\;\;\;\;\;\; {\rm if} \;\; \lambda^\prime \leq \lambda .\label{ortco}
\end{equation}
Combining \for{inco} and \for{ortco} entails
$\hat{D}_r\: p_\lambda =[\hat{D}_r]_{\lambda ,\lambda}\: p_{\lambda}$.
We now invoke Proposition~\ref{Deigv} to complete the proof.
\end{prf}

In Appendix~C it is shown that if a difference or differential operator
vanishes on ${\cal A}^W$, then all its coefficients must be zero.
Combining this result with Theorem~\ref{efunc} entails
the commutativity of the \ados{}:

\begin{thm}[commutativity]\label{commuting}\hfill

\noindent The operators $\hat{D}_1,\ldots ,\hat{D}_n$
mutually commute.
\end{thm}
\begin{prf}
The polynomials $\{ p_\lambda \}_{\lambda \in {\cal P}^+}$ form a basis of
${\cal A}^W$ consisting of joint eigenfunctions of
$\hat{D}_1, \ldots , \hat{D}_n$ (Theorem~\ref{efunc}).
Hence, it is clear that the \ados{} commute as operators on
${\cal A}^W$. In other words, the commutator of $\hat{D}_r$ and
$\hat{D}_{r^\prime}$
is a difference operator that vanishes on ${\cal A}^W$:
\begin{equation}
[\hat{D}_{r} ,\hat{D}_{r^\prime}] ( {\cal A}^W ) =0,\;\;\;\;\;\;\;\;\;
r, r^\prime  \in \{ 1,\ldots , n\} .
\end{equation}
It now follows from Proposition~\ref{vanish} in Appendix~C that the
coefficients of the commutator are identically zero.
\end{prf}

Consider the {\em real} algebra of difference operators
generated by $\hat{D}_1,\ldots ,\hat{D}_n$:
\begin{equation}
{\Bbb D} \equiv {\Bbb R}[\hat{D}_1,\ldots ,\hat{D}_n].
\end{equation}
It is clear that ${\Bbb D}$ is an abelian
algebra (Theorem~\ref{commuting}) and that the operators in ${\Bbb D}$
are simultaneously diagonalized by the Koornwinder polynomials
(Theorem~\ref{efunc}).

\begin{thm}[\mbox{${\Bbb D}\cong {\Bbb R}
[\ch\beta\theta_1,\ldots ,\ch\beta\theta_n ]^{S_n}$}]\label{hc} \hfill

\noindent For each symmetric function
$S(\theta )\in {\Bbb R}[\ch \beta \theta_1,\ldots , \ch \beta \theta_n]^{S_n}$
there exists a unique difference operator $\hat{D}\in {\Bbb D}$ such that
\begin{equation}
\hat{D}\: p_\lambda =S(\lambda +\rho )\: p_\lambda ,\;\;\;\;\;\;\;\;\;\;
\forall \lambda \in {\cal P}^+.\label{hciso}
\end{equation}
\end{thm}
\begin{prf}
$E_{r,n}(\theta )$ \for{ern}
is a linear combination of elementary symmetric functions:
\begin{equation}
E_{r,n}(\theta )=S_r(\ch\beta\theta_1,\ldots , \ch\beta\theta_n)+\;\;
            { l.d.},\;\;\;\;\;\;\;\;\;\;\;\;\;\; r=1,\ldots , n; \label{lincom}
\end{equation}
$l.d.$ stands for terms of lower degree in $\ch\beta\theta_j$, $j=1,\ldots ,
n$.
The elementary symmetric functions form a set of algebraically independent
generators of the symmetric algebra (this fact is the `fundamental
theorem on symmetric functions', see e.g. \cite{mac0}).
Hence, Eq. \for{lincom} implies that the same is true for the functions
$E_{r,n}(\theta )$: every element in ${\Bbb R}[\ch \beta \theta_1,\ldots ,
\ch \beta \theta_n]^{S_n}$ can be written uniquely as a polynomial
in $E_{r,n}(\theta )$, $r=1,\ldots , n$.

Now we use Theorem~\ref{efunc} to conclude that for every symmetric
function $S(\theta )$ there exists a difference operator
$\hat{D}\in {\Bbb D}$ such that Eq. \for{hciso} holds.
That such a difference operator $\hat{D}$ is unique follows from
Proposition~\ref{vanish} (Appendix~C).
\end{prf}

The symmetric functions
$S(\theta )\in {\Bbb R}[\ch \beta \theta_1,\ldots , \ch \beta \theta_n]^{S_n}$
separate the points of the wedge
\begin{equation}
 \{ \theta\in {\Bbb R}^n \; |\;
\theta_1\geq \theta_2\geq \cdots \geq \theta_n\geq 0 \} . \label{wedge}
\end{equation}
To see this, first notice that
$(-1)^{n-r}S_r(\ch\beta\theta_1,\ldots , \ch\beta\theta_n)$ is
the coefficient
of $\nu^{n-r}$ in the characteristic polynomial det$(T-\nu I)$ of the diagonal
matrix $T={\rm diag}(\ch\beta\theta_1,\ldots ,\ch\beta\theta_n)$.
Consequently, the values of $S_r(\ch\beta\theta_1,\ldots , \ch\beta\theta_n)$,
$r=1,\ldots , n$ determine $\theta$ in the wedge \for{wedge} uniquely.
This fact combined with Theorem~\ref{hc} can be used to prove the orthogonality
of the basis $\{ p_\lambda \} $:

\begin{cor}[orthogonality]\label{ort}
\begin{equation}
\langle p_{\lambda},p_{\lambda^\prime}\rangle_\Delta =0,\;\;\;\;\;\;\;\;\;\;\;
\forall \lambda ,\lambda^\prime \in {\cal P}^+,\; \lambda \neq \lambda^\prime .
\end{equation}
\end{cor}
\begin{prf}
Let $\lambda ,\lambda^\prime \in {\cal P}^+$, and $\lambda \neq
\lambda^\prime$.
Since the symmetric functions in
${\Bbb R}[\ch \beta \theta_1,\ldots , \ch \beta \theta_n]^{S_n}$ separate the
points of the wedge \for{wedge}, there exists an
$S(\theta )\in {\Bbb R}[\ch \beta \theta_1,\ldots , \ch \beta \theta_n]^{S_n}$
such that
\begin{equation}
S(\lambda +\rho )\neq S(\lambda^\prime +\rho ).
\end{equation}
Hence, by Theorem~\ref{hc} there exists a $\hat{D}\in {\Bbb D}$ for which
$p_\lambda$ and $p_{\lambda^\prime}$ are eigenfunctions
corresponding to different eigenvalues. But then the polynomials
$p_\lambda$ and $p_{\lambda^\prime}$ must be orthogonal
with respect to $\langle\cdot ,\cdot\rangle_\Delta$ because
$\hat{D}$ is symmetric, cf. Proposition~\ref{Dsym}.
\end{prf}

\vspace{1ex}
\noindent\note The orthogonality of the basis $\{ p_\lambda \}_{\lambda \in
{\cal P}^+}$
was already shown by Koornwinder \cite{koo3}.
His proof exploits the continuity of
$\langle p_{\lambda},p_{\lambda^\prime} \rangle_{\Delta}$
in the parameters.

\begin{cor}[self-adjointness]\label{sa}\hfill

\noindent Every difference operator $\hat{D}\in {\Bbb D}$ is essentially
self-adjoint on ${\cal A}^W\subset L^2_W({\Bbb T},\Delta dx)$.
\end{cor}
\begin{prf}
This an immediate consequence of the fact that every \ado{} in ${\Bbb D}$
acts as a real multiplication operator on the orthogonal basis
$\{ p_{\lambda}\}_{\lambda\in {\cal P}^+}$ of
$L^2_W({\Bbb T},\Delta dx)$.
\end{prf}

\vspace{1ex}
\noindent\remarks
{\em i.} Theorem~\ref{hc} states that the
assignments
\begin{equation}
\hat{D}_r \stackrel{HC}{\longmapsto} E_{r,n}(\theta ),\;\;\;\;\;\;
r=1,\ldots , n
\end{equation}
induce a Harish-Chandra-type algebra isomorphism
$HC:{\Bbb D}\stackrel{ \raisebox{-0.2ex}{$\sim$} }{\longrightarrow}
{\Bbb R}[\ch \beta \theta_1,\ldots , \ch \beta \theta_n]^{S_n}$.

{\em ii.} Recall that $p_\lambda$ is defined as $m_\lambda$ minus the
orthogonal
projection of $m_\lambda$ onto
span$\{ m_{\lambda^\prime} \}_{\lambda^\prime <\lambda }$.
Using the orthogonality of $\{ p_\lambda \}_{\lambda \in {\cal P}^+}$
(Corollary~\ref{ort}) this leads to the following recursion relation for
$p_\lambda$:
\begin{equation}
p_\lambda = m_\lambda -
\sum_{\lambda ^\prime\in {\cal P}^+,\: \lambda^\prime < \lambda}
\frac{\langle m_\lambda ,p_{\lambda^\prime}\rangle_\Delta}
     {\langle p_{\lambda^\prime},p_{\lambda^\prime}\rangle_\Delta}\;
p_{\lambda^\prime}.
\label{grs}
\end{equation}

{\em iii.} If the partial ordering~\ref{nat} is extended to a {\em linear}
ordering of the cone ${\cal P}^+$, then it is possible to orthogonalize the
basis $\{ m_\lambda \}$ by means of the Gram-Schmidt process.
By Corollary~\ref{ort}, the result
does not depend on the particular choice of the refinement of the ordering:
the resulting orthogonal basis coincides with $\{ p_\lambda \}$.
This amounts to a very restrictive property of the measure
$\Delta \, dx$.

{\em iv.} According to Theorem~\ref{commuting},
the Hamiltonians $\hat{D}_1,\ldots ,\hat{D}_n$
constitute a quantum integrable $n$-particle system
(cf. the note ending Section~\ref{ados}).

\setcounter{equation}{0}
\addtocontents{toc}{\protect\vspace{-2ex}}
\markright{4 $\beta \ra 0$: The Transition to $BC_n$-type Hypergeometric
\pdos{}}
\section{$\beta \ra 0$: The Transition to $BC_n$-type Hypergeometric
\pdos{}}\label{sect4}
By sending the step size $\beta$ of the differences to zero, our
\ados{} go over in commuting hypergeometric \pdos{} associated with the root
system $BC_n$. In this limit the eigenfunctions $\{ p_\lambda \}$ converge to
the $BC_n$-type Jacobi polynomials of Heckman and Opdam.

\noindent In this section we will make the dependence on $\beta$ explicit by
adding it as
a subscript on all objects of interest, e.g.: $\hat{D}_{r,\beta}$,
$\Delta_{\beta}$ and $p_{\lambda ,\beta}$.

\subsection{Eigenfunctions}\label{b=0eigf}
Consider the following weight function on the torus ${\Bbb T}$:
\begin{eqnarray}
\Delta_0(x) &\equiv&
\prod_{1\leq j<k\leq n}
 |\sin \alpha (x_j+x_k)\: \sin \alpha (x_j-x_k) |^{2g}
\nonumber \\
& & \prod_{1\leq j\leq n}
|\sin (\alpha x_j)|^{2\tilde{g}_0}\;
|\cos (\alpha x_j)|^{2\tilde{g}_1},
\;\;\;\;\;\;\;\;\;\;\;\;\; \alpha =1/2, \;\;
g, \tilde{g}_0, \tilde{g}_1 \geq 0 \label{deltanul}
\end{eqnarray}
and let $\langle \cdot ,\cdot \rangle_{\Delta_0}$ be the inner product
on $L^2_W({\Bbb T},\Delta_0\, dx)$ (cf. Eq. \for{ip}).
As before one introduces
$W$-invariant polynomials on ${\Bbb T}$ associated with the weight function
$\Delta_0$ (cf. Definition~\ref{koopol}); these are the $BC_n$-type Jacobi
polynomials of Refs.~\cite{heck1} and \cite{heck3}.
\begin{dfn}[$BC_n$-type Jacobi polynomials]\label{jacpol}\hfill

\noindent The {\em Jacobi polynomial}
$p_{\lambda ,0}\in {\cal A}^W_\lambda$
is defined by the conditions
\begin{equation}
p_{\lambda ,0}  = m_\lambda +
\sum_{\lambda^\prime \in {\cal P}^+,\, \lambda^\prime < \lambda}
c_{\lambda ,\lambda^\prime}\: m_{\lambda^\prime}  \label{jp1}
\end{equation}
and
\begin{equation}
\langle p_{\lambda ,0},m_{\lambda^\prime}\rangle_{\Delta_0} =0,
\;\;\;\;\;\; \forall \lambda^\prime \in {\cal P}^+,\;
\lambda^\prime <\lambda .\label{jp2}
\end{equation}
\end{dfn}

\vspace{1ex}
\noindent\remark Usually $\Delta_0$ is written in a slightly different form,
which
emphasizes the relation with the root system $BC_n$. Use
\begin{equation}
|\sin (\alpha x_j)|^{2\tilde{g}_0}\;
|\cos (\alpha x_j)|^{2\tilde{g}_1} \sim
|\sin (\alpha x_j)|^{2 k_0}\;
|\sin (2\alpha x_j)|^{2 k_1}
\end{equation}
with $\tilde{g}_0=k_0+k_1$ and $\tilde{g}_1=k_1$, to compare
\for{deltanul} with the usual expression for the $BC_n$-type weight function.

\vspace{1ex}
We need the following convergence result from \cite{koo1} to connect
$p_{\lambda ,0}$ with Koornwinder's polynomial
$p_{\lambda ,\beta}$, $\beta >0$:
\begin{lem}
Let $k_1, k_2 \in {\Bbb R}$ and $q\in ]0,1[$. Then
\begin{equation}
\lim_{q \uparrow 1}\: \frac{ (q^{k_1}z;q)_\infty }{ (q^{k_2}z;q)_\infty }
= (1-z)^{k_2-k_1},
\end{equation}
uniformly for $z$ in compacts of the punctured disc
$\{ z\in {\Bbb C}\; |\; |z|\leq 1,\; z\neq 1 \}$.
\end{lem}
\begin{prf}
See Proposition~{A.2} of Ref.~\cite[Appendix A]{koo1}.
\end{prf}
It is immediate from Eqs. \for{da+}, \for{db+f} and
the above lemma that
\begin{eqnarray}
\lim_{\beta \ra 0}\: d_{a,\beta}(z)&=& |2\sin (z/2)|^{2g}, \\
\lim_{\beta \ra 0}\: d_{b,\beta}(z)&=&
|2\sin (z/2)|^{2(g_0+g^\prime_0)}\; |2\cos (z/2)|^{2(g_1+g^\prime_1)},
\end{eqnarray}
uniformly for $z$ in compacts of $]\pi ,\pi [\setminus \{ 0\} $.

\noindent Consequently, for $\beta \ra 0$, $\Delta_\beta$ \for{densityb}
converges to a
weight function that is proportional to $\Delta_0$:
\begin{equation}
\lim_{\beta\ra 0}\: \Delta_\beta =
{\rm c}\: \Delta_0, \label{measlim}
\end{equation}
with
\begin{equation}
\tilde{g}_0=g_0+g^\prime_0 \: (\geq 0),\;\;\;\;\;\;\;\;\;\;\;\;
\tilde{g}_1=g_1+g^\prime_1 \: (\geq 0) \label{psub}
\end{equation}
and $c=4^{n(n-1)g} 4^{n(\tilde{g}_0+\tilde{g}_1)}$ ($>0$).

\begin{prp}\label{kootojac}
One has
\begin{equation}
\lim_{\beta \ra 0}\: p_{\lambda ,\beta}=p_{\lambda ,0},\;\;\;\;\;\;\;
\forall \lambda \in {\cal P}^+,\label{kojaeq}
\end{equation}
and
\begin{equation}
\langle p_{\lambda ,0},p_{\lambda^\prime ,0}\rangle_{\Delta_0}
=0,\;\;\;\;\;\;\;\;\;\;\;
\forall \lambda ,\lambda^\prime \in {\cal P}^+,\; \lambda \neq \lambda^\prime .
\end{equation}
\end{prp}
\begin{prf}
Using \for{grs}, \for{measlim} and an induction argument on $\lambda$, one
verifies that $\lim_{\beta \ra 0} p_{\lambda ,\beta}$ exists and satisfies
a recursion relation of the type \for{grs} (with
$\langle\cdot ,\cdot\rangle_{\Delta_\beta}\ra
 \langle\cdot ,\cdot\rangle_{\Delta_0}$).
Use Corollary~\ref{ort} and Eq.~\for{measlim} to conclude that the resulting
polynomials are orthogonal with respect to
$\langle \cdot ,\cdot \rangle_{\Delta_0}$.
But then the recursion relation for $\lim_{\beta \ra 0} p_{\lambda ,\beta}$
implies that it equals $m_\lambda$ minus the orthogonal projection
(with respect to $\langle \cdot ,\cdot \rangle_{\Delta_0}$) onto
span$\{ m_{\lambda^\prime} \}_{\lambda^\prime <\lambda}$.
Comparing with Definition~\ref{jacpol} one obtains \for{kojaeq}.
\end{prf}

\vspace{1ex}
\noindent\note The limit $\beta \ra 0$ corresponds to the limit $q\ra 1$,
cf. \for{psubsta}.
For Macdonald's polynomials the $q\ra 1$ limit to the Jacobi polynomials
of Heckman and Opdam was studied in \cite{mac1}, for arbitrary root systems.
The orthogonality of the Jacobi polynomials was proved in \cite{heck2}
(again for arbitrary root systems).

\subsection{Eigenvalues}
The purpose of this subsection is to investigate the behavior of the
eigenvalues of $\hat{D}_{r,\beta}$ (Proposition~\ref{Deigv})
for $\beta \ra 0$.
\begin{prp}\label{evbto0}
One has (for $r=1,\ldots ,n$)
\begin{eqnarray}
\lefteqn{\lim_{\beta \ra 0}\; \beta^{-2r}\,
          E_{r,n}(\ch\beta\theta_1,\ldots , \ch\beta\theta_n;
                   \ch\beta\rho_r,\ldots , \ch\beta\rho_n) }& & \nonumber \\
&=&  E_{r,n}(\theta_1^2/2,\ldots , \theta_n^2/2;
              \rho_r^2/2,\ldots , \rho_n^2/2)
= 2^{-r} E_{r,n}(\theta_1^2,\ldots , \theta_n^2;
        \rho_r^2,\ldots , \rho_n^2).\label{limE}
\end{eqnarray}
\end{prp}
\begin{prf}
Our proof of Eq. \for{limE} hinges on a recursion relation for $E_{r,n}$,
which has been relegated to Lemma~\ref{reclem2} of Appendix~B:
\begin{eqnarray}
\lefteqn{ E_{r,n}(\ch\beta\theta_1,\ldots ,
\ch\beta\theta_n;
                   \ch\beta\rho_r,\ldots , \ch\beta\rho_n)=} & & \nonumber \\
& &(\ch\beta \theta_n -\ch\beta\rho_n)\:
    E_{r-1,n-1}(\ch\beta\theta_1,\ldots , \ch\beta\theta_{n-1};
                   \ch\beta\rho_r,\ldots , \ch\beta\rho_n) \nonumber \\
& &+E_{r,n-1}(\ch\beta\theta_1,\ldots , \ch\beta\theta_{n-1};
                   \ch\beta\rho_r,\ldots , \ch\beta\rho_{n-1}) \;\;\;\;\;\;
				   1\leq r\leq n,\label{recrelbeta}
\end{eqnarray}
with the convention
\begin{equation}
E_{0,n}\equiv 1, \;\;\;\;\;\;\;\;\;\;
E_{r,n}\equiv 0 \;\; {\rm if}\;\; n<r .\label{inconv}
\end{equation}
We divide Eq. \for{recrelbeta} by $\beta^{2r}$ and then
use induction on $n$ to obtain
\begin{eqnarray}
\lefteqn{\lim_{\beta \ra 0} \beta^{-2r}
          E_{r,n}(\ch\beta\theta_1,\ldots , \ch\beta\theta_n;
                   \ch\beta\rho_r,\ldots , \ch\beta\rho_n) } \nonumber \\
&=& (\theta_n^2/2-\rho_n^2/2)\:
    E_{r-1,n-1}(\theta_1^2/2,\ldots , \theta_{n-1}^2/2;
              \rho_r^2/2,\ldots , \rho_n^2/2) \nonumber \\
& &+E_{r,n-1} (\theta_1^2/2,\ldots , \theta_{n-1}^2/2;
              \rho_r^2/2,\ldots , \rho_{n-1}^2/2). \label{recrel0}
\end{eqnarray}
Now we use again Lemma~\ref{reclem2} to conclude \for{limE} from relation
\for{recrel0}
(with convention \for{inconv}).
\end{prf}

\subsection{Operators}\label{oplim}
Expansion of $\hat{D}_{r,\beta}$ in $\beta$ yields a formal power
series of the form
\begin{equation}
\hat{D}_{r,\beta} =\sum_{m=0}^\infty
\hat{D}_{r}^{(m)} \beta^m . \label{opexp}
\end{equation}
The coefficients $\hat{D}_{r}^{(m)}$ are polynomials in the partials
$\hat{\theta}_j$, $j=1,\ldots , n$; this means that these coefficients are
\pdos{}.
We define the leading differential operator $\hat{D}_{r,0}$ of
$\hat{D}_{r,\beta}$ as the first nonzero coefficient in
expansion \for{opexp}:
\begin{dfn}[leading \pdo{}]\label{leadc} \hfill

\noindent Let
\begin{equation}
m_r \equiv
{\rm min} \{ m\in {\Bbb N}\; |\; \hat{D}_{r}^{(m)} \neq  0 \; \}
\end{equation}
(with $\hat{D}_{r}^{(m)}$ defined by expansion \for{opexp}).
Then,
\begin{equation}
\hat{D}_{r,0} \equiv \hat{D}_{r}^{(m_r)}
\end{equation}
is called the {\em leading \pdo{}} of $\hat{D}_{r,\beta}$.
\end{dfn}
The $BC_n$-type Jacobi polynomials are joint eigenfunctions of
$\hat{D}_{1,0},\ldots  ,\hat{D}_{n,0}$:
\begin{thm}\label{ev0}
One has $m_r=2r$ and
\begin{equation}
\hat{D}_{r,0}\: p_{\lambda ,0}=
E_{r,n}\left( (\lambda_1+\rho_1)^2,\ldots , (\lambda_n+\rho_n)^2;
        \rho_r^2,\ldots , \rho_n^2\right) \: p_{\lambda
,0},\;\;\;\;\;\;\;\;\;\;\;
		\forall \lambda \in {\cal P}^+,\label{eveq0}
\end{equation}
with $r=1,\ldots ,n$ and $\rho$ as in \for{roo}.
\end{thm}
\begin{prf}
Consider the eigenvalue equation \for{Dreveq}:
\begin{equation}
\hat{D}_{r,\beta}\:(p_{\lambda ,\beta})=2^r\,
E_{r,n}(\ch\beta (\lambda_1+\rho_1),\ldots , \ch\beta (\lambda_n+\rho_n);
        \ch\beta\rho_r,\ldots , \ch\beta\rho_n)\: p_{\lambda
,\beta}.\label{eveq}
\end{equation}
First, apply Taylor's theorem to the l.h.s. of Eq. \for{eveq} and
make use of Definition~\ref{leadc} and Limit~\for{kojaeq} to conclude that
\begin{equation}
\hat{D}_{r,\beta}\, p_{\lambda ,\beta}=
\hat{D}_{r,0}\, p_{\lambda ,0}\; \beta^{m_r}\;\;\;
+o(\beta^{m_r}).\label{lhseq}
\end{equation}
Next, use Proposition~\ref{evbto0} and Limit~\for{kojaeq} to derive
the asymptotic behavior for $\beta \ra 0$ of the r.h.s. of \for{eveq}:
\begin{eqnarray}
\lefteqn{ 2^r\,
 E_{r,n}(\ch\beta (\lambda_1+\rho_1),\ldots ,\ch\beta (\lambda_n+\rho_n);
 \ch\beta\rho_r,\ldots , \ch\beta\rho_n)\: p_{\lambda ,\beta} =} & & \nonumber
\\
& & E_{r,n}\left( (\lambda_1+\rho_1)^2,\ldots ,(\lambda_n+\rho_n)^2;
          \rho_r^2,\ldots , \rho_n^2 \right)\: p_{\lambda ,0}\; \beta^{2r}
\;\;\; +o(\beta ^{2r}). \label{rhseq}
\end{eqnarray}
By Definition~\ref{leadc} and Proposition~\ref{vanish} it is possible to pick a
$\lambda \in {\cal P}^+$
such that $\hat{D}_{r,0}\, p_{\lambda ,0}\neq 0$, so $m_r\geq 2r$.
It is not difficult to see that there also exist $\lambda \in {\cal P}^+$
such that
\begin{displaymath}
E_{r,n}\left( (\lambda_1+\rho_1)^2,\ldots ,(\lambda_n+\rho_n)^2;
          \rho_r^2,\ldots , \rho_n^2 \right) \neq 0,
\end{displaymath}
 so $m_r \leq 2r$.
This entails $m_r=2r$ and \for{eveq0}.
\end{prf}
\begin{cor}
\begin{equation}
\hat{D}_{r,0} =\lim_{\beta \ra 0}\; \beta^{-2r}\:
       \hat{D}_{r,\beta}.\label{oplimit}
\end{equation}
\end{cor}
Explicit computation of \for{oplimit} for $r=1$ yields
(use \for{D1}, \for{parres})
\begin{eqnarray}
\hat{D}_{1,0}
&=& \sum_{1\leq j\leq n} \hat{\theta}_j^2
-2i\alpha g \sum_{1\leq j<k\leq n} \left\{
\cot\alpha (x_j+x_k) (\hat{\theta}_j+\hat{\theta}_k)+
\cot\alpha (x_j-x_k) (\hat{\theta}_j-\hat{\theta}_k)\right\} \nonumber \\
& &\;\;\;\;\;\;\;\;\;\;\;\;\;\;\; -2i\alpha \sum_{1\leq j\leq n}
\left\{  \tilde{g}_0\cot (\alpha x_j) -\tilde{g}_1\tan (\alpha x_j) \right\}
\hat{\theta}_j \label{D10}
\end{eqnarray}
(with $\alpha =1/2$ and $\tilde{g}_0$, $\tilde{g}_1$ as in Eq. \for{psub}).

An immediate consequence of Theorem~\ref{commuting} and Limit~\for{oplimit} is
the commutativity of $\hat{D}_{r,0}$, $r=1,\ldots ,n$.
\begin{cor}\label{compdo}
The differential operators $\hat{D}_{1,0},\ldots ,\hat{D}_{n,0}$,
mutually commute.
\end{cor}
Let ${\Bbb D}_0\equiv {\Bbb R}[\hat{D}_{1,0},\ldots ,\hat{D}_{n,0}]$.
The algebra ${\Bbb D}_0$ consists of commuting \pdos{} which are simultaneously
diagonalized by the $BC_n$-type Jacobi polynomials $p_{\lambda ,0}$.

\begin{thm}[\mbox{${\Bbb D}_0 \cong
{\Bbb R}[\theta_1^2,\ldots , \theta_n^2]^{S_n}$}]\label{hc0}\hfill

\noindent For every symmetric polynomial $S(\theta )\in
{\Bbb R}[\theta_1^2,\ldots , \theta_n^2]^{S_n}$
there exists a unique differential operator $\hat{D}_0\in {\Bbb D}_0$ such
that
\begin{equation}
\hat{D}_0\: p_{\lambda ,0}=S(\lambda +\rho )\: p_{\lambda ,0},
\;\;\;\;\;\;\;\;\;\; \forall \lambda \in {\cal P}^+.
\end{equation}
\end{thm}
\begin{cor}[self-adjointness]\label{sa0}\hfill

\noindent The \pdos{} $\hat{D}_0\in {\Bbb D}_0$ are essentially
self-adjoint on ${\cal A}^W\subset L^2_W({\Bbb T},\Delta_0dx)$.
\end{cor}
The proofs of Theorem~\ref{hc0} and Corollary~\ref{sa0} are virtually the same
as for
Theorem~\ref{hc} and Corollary~\ref{sa}, respectively.

\vspace{1ex}
\remarks {\em i.} The operator $\hat{D}_{1,0}$ \for{D10} coincides with the
simplest (i.e. lowest order) hypergeometric \pdo{} associated with the root
system $BC_n$. In order to compare \for{D10} with the usual notation which
emphasizes the r\^ole of the root system, one should eliminate
$\tan (\alpha x_j)$ from \for{D10} by means of the relation
$\cot (2\alpha x_j)= ( \cot (\alpha x_j)-\tan (\alpha x_j) )/2$
(cf. the remark under Definition~\ref{jacpol}).

{\em ii.} The existence of an abelian algebra of \pdos{} containing
$\hat{D}_{1,0}$, which is isomorphic to
${\Bbb R}[\theta_1^2,\ldots ,\theta_n^2]^{S_n}$ via a
Harish-Chandra-type isomorphism (cf. Theorem~\ref{hc0}), was already shown by
Heckman and Opdam \cite{heck2,opd2}.
In fact, their result is more general since they consider arbitrary root
systems.

{\em iii.} For an arbitrary root system $R$, the quantum Hamiltonian $\hat{H}$
of
the corresponding generalized Calogero-Sutherland system is related to the
second order hypergeometric \pdo{} via a similarity transformation \cite{op2}.
{}From a mathematical point of view, this just amounts to the conjugation with
$\Delta ^{1/2}$, which transforms between Lebesgue measure and
Plancherel measure with weight function $\Delta$.
In our case (i.e. for $R=BC_n$), one has
\pagebreak[2]
\begin{eqnarray}
\hat{H}&=&\Delta_0^{1/2} \hat{D}_{1,0}\, \Delta_0^{-1/2}\; +\, E_0  \nonumber
\\
&=& \sum_{1\leq j\leq n} \hat{\theta}_j^2 \; + \;
2g(g-1)\, \alpha^2  \sum_{1\leq j< k\leq n} \left\{
\sin^{-2} \alpha (x_j+x_k)\; +\; \sin^{-2} \alpha (x_j-x_k) \right\}  \nonumber
\\
& & +\,\alpha^2 \sum_{1\leq j\leq n} \left\{
\tilde{g}_0(\tilde{g}_0-1)\: \sin^{-2} (\alpha x_j) \; + \;
\tilde{g}_1(\tilde{g}_1-1)\: \cos^{-2} (\alpha x_j) \right\}
\end{eqnarray}
with $E_0=4\alpha^2 (\rho ,\rho )$ and
$\rho_j= (n-j)g+ (\tilde{g}_0+\tilde{g}_1)/2$, cf. \for{roo}.

\noindent Corollary~\ref{compdo} can be interpreted as the quantum
integrability
of the $BC_n$-type Calogero-Sutherland system.
For arbitrary root systems integrability follows from \cite{heck2,opd2}
(cf. Remark {\em ii}).

\setcounter{equation}{0}
\addtocontents{toc}{\protect\vspace*{-1ex}}
\markright{5 Special Cases Related to Classical Root Systems}
\section{Special Cases Related to Classical Root Systems}\label{sect5}
By limit transitions and/or specialization
of the parameters, the operators $\hat{D}_1,\ldots ,\hat{D}_n$ reduce to
commuting \ados{}, which are simultaneously diagonalized by Macdonald's
polynomials. Such difference operators are obtained
for all Macdonald families
associated with (admissible pairs of) the classical root systems:
$A_{n-1}$, $B_n$, $C_n$, $D_n$ and $BC_n$.

\vspace{1ex}
\noindent\note Most results in this section have an obvious counterpart for
$\beta = 0$ (which amounts to $q= 1$).

\subsection{Preliminaries}
First, we outline very briefly some of the main points of the construction
presented by Macdonald \cite{mac1}.
A more detailed summary of his results can be found in \cite{mac2} and
\cite{koo3}. (For our purposes, especially the second summary is useful).
Here, we only want to introduce some terminology
which facilitates clarifying the connection between the preceding
sections and Ref.~\cite{mac1}. For general information on root systems
the reader is referred to e.g. \cite{bou,ser}.
Although most of the remaining part of the paper
should be accessible without a detailed knowledge of root systems,
a glance at the `planches' in Bourbaki \cite{bou}
might be of some help.

Ref.~\cite{mac1} uses the concept of admissible pairs of root systems.
The pair
$(R,S)$ is admissible if $R$ and $S$ are root systems
(assumed irreducible) such that
$S\subset R$ is reduced and generates the same Weyl group as $R$.
Let $V$ be the real vector space spanned by $R$ and
consider the torus ${\Bbb T}_R \equiv V /(2\pi {\Bbb Z} R^\vee )$.
Let ${\cal P}^+_R$ be the dominant cone of the weight lattice of $R$
(which equals the character lattice of ${\Bbb T}_R$) and let
${\cal A}^W_R$ denote the algebra of W-invariant (trigonometric) polynomials on
${\Bbb T}_R$ (this algebra is isomorphic to the W-invariant part
of the group algebra over the weight lattice).
To every admissible pair $(R,S)$ Macdonald associates a weight function
$\Delta_{(R,S)}$ on ${\Bbb T}_R$ and finds a corresponding orthogonal basis
$\{ p_{\lambda , (R,S)} \}_{\lambda \in {\cal P}_R^+}$ of ${\cal A}^W_R$.
Furthermore, he introduces difference operators $D_\sigma$ associated with the
so-called (quasi-)minuscule weights $\sigma$ of $S^\vee$. These operators
are diagonalized by the basis of Macdonald polynomials
$\{ p_{\lambda , (R,S)} \}_{\lambda \in {\cal P}_R^+}$.

We will show that additional \ados{} for Macdonald's polynomials
arise as special cases of
$\hat{D}_1,\ldots ,\hat{D}_n$ \for{Dr1}.
This leads to difference operators associated with the fundamental
weights of $S^\vee$, for every admissible pair consisting
of classical root systems.
These \ados{} generate an abelian algebra ${\Bbb D}_{(R,S)}$ of
difference operators, which are simultaneously diagonalized by the basis
$\{ p_{\lambda , (R,S)} \}_{\lambda \in {\cal P}_R^+}$.
The algebra ${\Bbb D}_{(R,S)}$ is isomorphic to ${\Bbb R}[S^\vee ]^W$
(the W-invariant part of the real group algebra over ${\cal P}_S^\vee$).

\vspace{1ex}
\noindent\remark If all roots in $R$ have the same length
(this is the case for $R=A_n$ or $D_n$), then $S=R$ and
there exists only one admissible pair.
If there are roots with different lengths, then there are several possibilities
for the pair $(R,S)$. In case of the classical series there are six such
possibilities; these correspond to $R=B_n$, $C_n$ or $BC_n$ and $S=B_n$ or
$C_n$.

\subsection{The Root System $A_{n-1}$}\label{mdonaldAn}
Let $\hat{D}_{r,lead}$ consist of those terms in $\hat{D}_r$
that are of highest order
in the exponentials $\exp (-\beta\hat{\theta}_j)$, $j=1,\ldots , n$
(cf. Eq.~\for{Dr3}):
\begin{equation}
\hat{D}_{r,lead} = \sum_{\stackrel{J\subset \{ 1,\ldots , n\}}{ |J|=r }}
V_{J;J^c}\: e^{-\beta \hat{\theta}_J},\;\;\;\;\;\;\;\;\;\;
r=1,\ldots , n.
\end{equation}
One picks up these leading terms via the following limit:
\begin{equation}
\hat{D}_{r,lead}= \lim_{R \ra \infty}\; e^{-r\beta R}\;
\left( \Lambda_R^{-1}\, \hat{D}_{r}\, \Lambda_R \right) ,\label{leadtado}
\end{equation}
with
\begin{equation}
\Lambda_R \equiv e^{-iR(x_1+\cdots +x_n)}.
\end{equation}
Let
\begin{equation}
\tilde{\Delta}_+ \equiv
\prod_{1\leq j<k\leq n} \left( e^{g(x_j+x_k)/2i}\: d_a^+(x_j+x_k) \right)\;
\prod_{1\leq j\leq n} \left( e^{(g_0+g_1+g_0^\prime +g_1^\prime )x_j/2i}\:
d_b^+(x_j) \right) .
\end{equation}
Conjugation of $\hat{D}_{r,lead}$ with $\tilde{\Delta}_+$ results
in $A_{n-1}$-type difference operators (use \for{d1}, \for{d2}):
\begin{eqnarray}
\hat{D}^\prime _{r}&=& \tilde{\Delta}_+\, \hat{D}_{r,lead}\,
\tilde{\Delta}_+^{-1} \label{conju} \\
&=& \sum_{\stackrel{J\subset \{ 1,\ldots , n\}}{ |J|=r }}
\left( \prod_{\stackrel{j\in J}{k\in J^c}} v_a(x_j-x_k) \right)
 e^{-\beta \hat{\theta}_J},\;\;\;\;\;\;\;\;\;\;
r=1,\ldots , n.\label{Drprime}
\end{eqnarray}
It is clear from Theorem~\ref{commuting} and Eqs.
\for{leadtado} and \for{conju} that the operators
$\hat{D}^\prime_{1},\ldots ,\hat{D}^\prime_{n}$ commute.
For $r=n$, \for{Drprime} reduces to an operator that only generates
a translation of the coordinates:
$\hat{D}^\prime_{n}=\exp\left( -\beta (\hat{\theta}_1+\cdots
                +\hat{\theta}_n) \right) $.
For $r<n$, $\hat{D}^\prime_{r}$ decomposes in two
commuting parts:
\begin{equation}
\hat{D}^\prime_{r}=  e^{-\beta r(\hat{\theta}_1+\cdots +\hat{\theta}_n)/n}
                     \; \hat{D}_{r, A_{n-1}},
\;\;\;\;\;\;\;\;\;\;\; 1\leq r\leq n-1.\label{DrAn}
\end{equation}
The first part causes a translation
$x\ra x+i\beta (r/n)\: (e_1+\cdots +e_n)$;
the second part coincides (up to a multiplicative constant)
with Macdonald's difference operator $E_{\omega_r}$:
\begin{equation}
E_{\omega_r}= c_r\: \hat{D}_{r, A_{n-1}},
\;\;\;\;\; \;\;\; \;\;\;\;\; c_r\equiv e^{-\beta gr(n-r)/2 },
\;\;\;\;\;\;\;\;\;\;\;\;\;\;\; 1\leq r\leq n-1. \label{mdA}
\end{equation}
The operator $E_{\omega_r}$ is associated with the $r$th fundamental
weight $\omega _r$ of the root system $A_{n-1}$.
The parameters in \cite{mac1} are related to ours via Eq. \for{psubsta}.

\vspace{1ex}
\noindent\notes {\em i.} For $r=1$, Eq.~\for{Drprime} reduces to
\begin{equation}
\hat{D}^\prime_{1} = \sum_{1\leq j\leq n}
\left( \prod_{k\neq j} v_a(x_j-x_k) \right) \; e^{-\beta\hat{\theta}_j}.
\end{equation}

{\em ii.} After transformation to Lebesgue measure, the operator
$\hat{D}^\prime_{r}$
goes over in the $r$\-th quantum integral $\hat{S}_r$ \cite[(Eq. (2.3))]{rui2}
of the relativistic Calogero-Moser system with trigonometric coefficients.
More precisely, let
\begin{equation}
\Delta_{A_{n-1}}(x) \equiv \prod_{1\leq j<k\leq n} d_a(x_j-x_k),
\end{equation}
then
\begin{eqnarray}
\hat{S}_r &=& \Delta_{A_{n-1}}^{1/2}\: \hat{D}^\prime _{r}\:
\Delta_{A_{n-1}}^{-1/2} \nonumber \\
&=& \sum_{\stackrel{J\subset \{ 1,\ldots , n\}}{ |J|=r }}
\left( \prod_{\stackrel{j\in J}{k\in J^c}} v_a(x_j-x_k) \right)^{1/2}
e^{-\beta \hat{\theta}_J}
\left( \prod_{\stackrel{j\in J}{k\in J^c}} v_a(x_k-x_j) \right)^{1/2} .
\end{eqnarray}
This relation between the $n$-particle relativistic $CM$ system introduced by
Ruijsenaars and Macdonald's
difference operators for the root system $A_{n-1}$ was first observed by
Koornwinder \cite{koo0}. It generalizes the relation between the $n$-particle
Calogero-Sutherland system and the hypergeometric \pdos{} associated with
$R=A_{n-1}$ \cite{op2}
(cf. Remark {\em iii} at the end of Section~\ref{oplim}).

{\em iii.} The operators $\hat{D}_{r,A_{n-1}}$ act as a real multiplication
on functions that depend only on $x_1+\cdots +x_n$.
The transition
$\hat{D}^\prime_{r}\longrightarrow \hat{D}_{r,A_{n-1}}$
can be interpreted physically as restricting attention to the motion in the
center of mass hyperplane $x_1+\cdots +x_n=0$.

\vspace{1ex}
We show next that the joint eigenfunctions of $\hat{D}_{r,A_{n-1}}$,
i.e. Macdonald's $A_{n-1}$-type polynomials, can be obtained from
Koornwinder's polynomials by a certain limit transition.
Let $m_{\lambda ,lead}$ be the sum of terms in $m_\lambda$ \for{monomial}
which are of the
highest degree in $\exp (ix_j)$, $j=1,\ldots , n$:
\begin{equation}
m_{\lambda ,lead} = \sum_{\lambda^\prime \in S_n\lambda} e^{\lambda^\prime},
\;\;\;\;\;\;\;\;\; \lambda \in {\cal P}^+.
\end{equation}
Recall that according to Definition~\ref{koopol}, $p_\lambda$
is a linear combination of monomials of the form
\begin{equation}\label{pl}
p_\lambda =
\sum_{\lambda^\prime \in {\cal P}^+,\, \lambda^\prime \leq \lambda}\;
c_{\lambda ,\lambda^\prime}\;
m_{\lambda^\prime}\; ,\;\;\;\;\;\;\;\;\;\;\;\;\lambda\in {\cal P}^+,
\end{equation}
with $c_{\lambda ,\lambda^\prime}$ certain complex coefficients
(which depend only on the parameters) such that \for{kp2} holds
and $c_{\lambda ,\lambda}=1$.
We set
\begin{equation}\label{plead}
p_{\lambda ,lead}\equiv
\sum_{\stackrel{\lambda^\prime\in {\cal P}^+,\,\lambda^\prime \leq \lambda}
               {|\lambda^\prime |=|\lambda |}}\;
c_{\lambda ,\lambda^\prime}\;
m_{\lambda^\prime ,lead} \; ,\;\;\;\;\;\;\;\;\;\lambda\in {\cal P}^+.
\end{equation}
Let
\begin{equation}
\omega \equiv e_1+\cdots +e_n. \label{omga}
\end{equation}
{}From the asymptotics
\begin{equation}
m_\lambda \left( x-iR \omega \right) \sim
m_{\lambda ,lead}(x)\: e^{R |\lambda | },\;\;\;\;\;\; R \ra \infty,
\end{equation}
and Eqs. \for{pl} and \for{plead}, one derives
\begin{equation}
p_{\lambda ,lead}(x) =\lim_{R\ra\infty}\; e^{-R |\lambda | }\:
p_{\lambda}\left( x-iR\omega \right) .\label{kootomac}
\end{equation}
The polynomial $p_{\lambda ,lead}$ is homogeneous
of degree $|\lambda |$ in the exponentials $\exp (ix_j)$.
Consequently, a translation causes
an automorphic phase factor: $x\ra x+R\omega$ $\Rightarrow$
$p_{\lambda ,lead} \ra \exp (i R |\lambda |)\; p_{\lambda ,lead}$.
By multiplying
$p_{\lambda ,lead}$ with an appropriate exponential function one ends up
with a basis of translation-invariant functions:
\begin{equation}
p_{\lambda ,A_{n-1}}(x) \equiv
e^{-i|\lambda |\: (x_1+\cdots +x_n)/n}\: p_{\lambda ,lead}(x) =
\sum_{\stackrel{\lambda^\prime\in {\cal P}^+,\,\lambda^\prime \leq \lambda}
               {|\lambda^\prime |=|\lambda |}}\;
c_{\lambda ,\lambda^\prime}\;
m_{\lambda^\prime ,A_{n-1}}(x),\label{macdAnpols}
\end{equation}
with
\begin{equation}
m_{\lambda^\prime ,A_{n-1}}(x)\equiv
e^{-i|\lambda^\prime |\: (x_1+\cdots +x_n)/n}\: m_{\lambda^\prime ,lead}(x).
\end{equation}
As the notation suggests, it will turn out that the functions
$\{ p_{\lambda ,A_{n-1}}\}$ \for{macdAnpols} coincide with the
Macdonald polynomials associated with the root system $A_{n-1}$.

First we need a lemma. It says that the spectrum of
$\hat{D}_1$ \for{D1} is monotonic in $\lambda\in {\cal P}^+$.
(Recall that the eigenvalues of $\hat{D}_1$ are given by \for{ev} with $r=1$;
see also Eq.~\for{D1eigen} below).
\begin{lem}\label{nondeg}
Let
\begin{equation}
F_\beta (\theta )\equiv \sum_{1\leq j\leq n} \ch\beta \theta_j.
\end{equation}
Then, for all $\lambda ,\lambda^\prime \in {\cal P}^+$:
\begin{equation}
\lambda  > \lambda ^\prime \Longrightarrow
F_\beta (\lambda +\rho ) > F_\beta (\lambda^\prime +\rho )\label{monotone}
\end{equation}
(with $\rho$ given by \for{roo}).
\end{lem}
\begin{prf}
It is clear from Definition~\ref{nat} that $\lambda > \lambda ^\prime$ iff
\begin{equation}
\lambda =\lambda^\prime +\sum_{1\leq j\leq n-1}
a_j\, (e_j-e_{j+1}) \; +\; a_n\,e_n, \;\;\;\;\;\;\;\;\;\;\;\;\;\;
a_j \in {\Bbb N}
\end{equation}
with at least one of the $a_j$'s positive.

Obviously, it suffices to verify \for{monotone} for the special case that
only one of the $a_j$ is positive. Now for $j=n$ this is immediate,
while for $j<n$ this follows because $\ch\beta (\cdot )$ is a convex
function:
\begin{equation}
\ch\beta (x +a) +\ch\beta (y -a) > \ch\beta x +\ch\beta y
\end{equation}
if $x\geq y$ and $a>0$.
\end{prf}

\begin{prp}\label{mcApol}
Let $\hat{D}_{r,A_{n-1}}$ be determined by \for{Drprime}-\for{DrAn}
and let
$p_{\lambda ,A_{n-1}}$ be defined by \for{macdAnpols}. Then
\begin{equation}
\hat{D}_{r,A_{n-1}}\: p_{\lambda , A_{n-1}}\; =\;
E_{r,A_{n-1}}(\lambda +\rho^\prime )\: p_{\lambda ,A_{n-1}},\label{eveqAn}
\end{equation}
with
\begin{eqnarray}
E_{r,A_{n-1}}(\theta )&\equiv &
e^{\beta r (\theta_1+\cdots +\theta_n)/n}\;
S_r\left( e^{-\beta \theta_1},\ldots ,e^{-\beta \theta_n}\right) , \\
\rho_j^\prime &\equiv &g(n+1-2j)/2, \;\;\;\;\;\;\;\; 1\leq j\leq
n.\label{rooprime}
\end{eqnarray}
($S_r$ denotes the rth elementary symmetric function (Definition~\ref{elsy})).
\end{prp}
\begin{prf}
The operator $\hat{D}^\prime_r$ \for{Drprime} is invariant both
under permutations
of $x_j$ and under translations of the form $x\ra x+R\omega$
(with $\omega$ as in \for{omga}).
We use this and the asymptotics of
$(\hat{D}_r^\prime m_{\lambda ,lead})(-iRy)$ for $R\ra \infty$
(with $y$ such that \for{chamber} holds) to derive
(cf. Proposition~\ref{inv}, Proposition~\ref{Dtri}, and their proofs)
\begin{equation}\label{trlead}
\hat{D}_r^\prime m_{\lambda ,lead} =
\sum_{\stackrel{\lambda^\prime \in {\cal P}^+,\, \lambda^\prime \leq \lambda}
                {|\lambda ^\prime | =|\lambda |}}
[\hat{D}_r^\prime ]_{\lambda ,\lambda^\prime}\; m_{\lambda^\prime ,lead}\; ,
\end{equation}
with
\begin{equation}\label{evlead}
[\hat{D}_r^\prime ]_{\lambda ,\lambda} =
S_r(e^{-\beta (\lambda_1 +\rho_1^\prime)},\ldots ,
                     e^{-\beta (\lambda_n +\rho_n^\prime )}).
\end{equation}
(The poles at $x_j=x_k$, $j\neq k$, cancel because of the permutation symmetry;
the condition $|\lambda^\prime |=|\lambda |$ in sum the \for{trlead} stems from
the translational invariance of $\hat{D}_r^\prime$).

We will now show that $p_{\lambda ,lead}$ is an eigenfunction of
$\hat{D}_r^\prime$ (with eigenvalue \for{evlead}).
Consider the eigenvalue equation \for{Dreveq} for $r=1$:
\begin{equation}
(\hat{D}_{1}\, p_{\lambda})(x)=
2\left(\sum_{1\leq j\leq n}[\ch\beta (\lambda_j+\rho_j)-\ch\beta\rho_j]\right)
\: p_{\lambda}(x),\;\;\;\;\;\;\;\;\;\; \lambda \in {\cal P}^+,\label{D1eigen}
\end{equation}
with $\hat{D}_1$ given by \for{D1}.
Substitute
\begin{equation}
x\ra x -iR\omega
\end{equation}
and divide both sides of the equation by
$\exp (R |\lambda |)$; Sending $R\ra\infty$ entails (use
\for{asva},\for{asvb} and \for{kootomac})
\begin{eqnarray}
\lefteqn{
\left( c_1\, \hat{D}^\prime_{1} +
c_1^{-1}\, (\hat{D}^\prime_{n})^{-1}\hat{D}^\prime_{n-1}\; -c_2\right) \:
p_{\lambda ,lead} = } & & \nonumber \hfill \\
& & 2\left(\sum_{1\leq j\leq n}[\ch\beta
(\lambda_j+\rho_j)-\ch\beta\rho_j]\right)\:
p_{\lambda ,lead},\label{limeig}
\end{eqnarray}
with
\begin{eqnarray}
c_1 &=& e^{-\beta g(n-1)/2} e^{-\beta (g_0+g_1+g_0^\prime +g_1^\prime )/2},\\
c_2 &=& c_1\sum_{1\leq j\leq n} \prod_{k\neq j}v_a(x_j-x_k)+
   c_1^{-1}\sum_{1\leq j\leq n} \prod_{k\neq j}v_a(x_k-x_j) \label{regAn} \\
&\stackrel{*}{=}& 2\sum_{1\leq j\leq n}\ch\beta\rho_j .
\end{eqnarray}
(To verify equality *, first check that both parts of \for{regAn} are regular
and
bounded in $x_j$; consequently, these parts are constants because of
Liouville's theorem.
One obtains the value of these constants
by putting $x=iRy$ with $y$ such that \for{chamber} holds, and
then sending $R\ra \infty$).

It is clear from the commutativity of
$\hat{D}_1^\prime ,\ldots ,\hat{D}_n^\prime$ that $\hat{D}_r^\prime$ commutes
with the operator on the l.h.s. of \for{limeig}.
Therefore, $\hat{D}_r^\prime p_{\lambda ,lead}$ is an eigenfunction
of the latter operator corresponding to the same eigenvalue as
$p_{\lambda ,lead}$. We know from \for{plead} and \for{trlead} that
$\hat{D}_r^\prime p_{\lambda ,lead}$ lies  in $ {\rm span}
\{ \; p_{\lambda^\prime ,lead} \; |\; \lambda^\prime \leq \lambda,\;
|\lambda^\prime |=|\lambda |\; \}$.
Now we use Lemma~\ref{nondeg} to conclude that
$p_{\lambda ,lead}$ must be an eigenfunction of
$\hat{D}_r^\prime $.
The corresponding eigenvalue follows from \for{trlead}, \for{evlead}.

One obtains the  expressions \for{eveqAn}-\for{rooprime} by restricting to the
hyperplane $x_1+\cdots +x_n=0$.
\end{prf}
Notice that $m_{\lambda^\prime ,A_{n-1}}=m_{\lambda ,A_{n-1}}$ iff
$\lambda^\prime -\lambda \in {\Bbb Z}(e_1+\cdots +e_n)$. Thus,
the monomials $m_{\lambda ,A_{n-1}}$ can be relabeled by the projection
of ${\cal P}^+$ \for{domcone} onto the hyperplane $x_1+\cdots +x_n=0$.
This projection of ${\cal P}^+$ coincides with the cone
${\cal P}^+_{A_{n-1}}$ of dominant $A_{n-1}$ weight vectors.
The polynomials $p_{\lambda ,A_{n-1}}$ can also be relabeled by
${\cal P}^+_{A_{n-1}}$, since
$p_{\lambda^\prime ,A_{n-1}}=p_{\lambda ,A_{n-1}}$ iff
$\lambda^\prime -\lambda \in {\Bbb Z}(e_1+\cdots +e_n)$.
(This follows from expansion \for{macdAnpols} and the eigenvalue equations
\for{eveqAn} for $r=1,\ldots ,n-1$).
Since the operators $\hat{D}_{r, A_{n-1}}$ coincide with Macdonald's
$A_{n-1}$ difference operators up to a constant,
we end up with the following corollary:
\begin{cor}
The function $p_{\lambda ,A_{n-1}}$ coincides with
the $A_{n-1}$-type Macdonald polynomial corresponding to the
weight vector $\lambda -|\lambda |(e_1+\cdots +e_n)/n \in {\cal P}^+_{A_{n-1}}$
(with the parameters as in \for{psubsta}).
\end{cor}

\vspace{1ex}
\noindent\notes {\em i.}
The transition $p_{\lambda ,lead}\ra p_{\lambda ,A_{n-1}}$
amounts to ignoring the linear motion of the center of mass.

{\em ii.} For $\beta \ra 0$ (i.e. $q=\exp (-\beta )\ra 1$), the polynomials
 $p_{\lambda ,A_{n-1}}$ converge to the
Jacobi polynomials associated with $A_{n-1}$ \cite{mac1}.
It is clear that the $\beta = 0$ version of the formulas \for{plead},
\for{macdAnpols} relates
the Jacobi polynomials associated with $BC_n$ and $A_{n-1}$:
\begin{equation}
p_{\lambda ,A_{n-1}}= e^{-i|\lambda | (x_1+\cdots +x_n)/n}\;
\lim_{R\ra \infty}\;\;\;
e^{-R|\lambda |}\: p_{\lambda ,BC_n} (x-iR(e_1+\cdots +e_n)).
\end{equation}

{\em iii.} Recently, a completely different limit taking the Jacobi polynomials
associated
with $BC_n$ to those associated with $A_{n-1}$ has been found \cite{bk};
this limit has not been dealt with for $q\neq 1$.

\subsection{The Root Systems $B_n$, $C_n$ and $BC_n$.}
In order to compare our results with \cite{mac1},
it is convenient to carry out a reparametrization:
\begin{equation}
\begin{array}{lcl}
\mu_0 \ra \nu_1 +\nu_2, &\;  &\mu_0^\prime \ra \nu_1^\prime +\nu_2^\prime ,\\
\mu_1 \ra \nu_2,      &\;  &\mu_1^\prime \ra \nu_2^\prime
\end{array}
\end{equation}
with (cf. \for{parres})
\begin{equation}
\nu_\delta \equiv i\beta k_\delta ,\;\;\;\;\;
\nu_\delta^\prime \equiv i\beta k_\delta^\prime ,\;\;\;\;\;\;
k_\delta ,k_\delta^\prime \geq 0, \;\;\;\;\;\;
\delta =1,2. \label{parresred}
\end{equation}
With these new parameters we rewrite $v_b(z)$ \for{vb} and $d_b^+(z)$
\for{db+} (recall also \for{db+f}):
\begin{eqnarray}
\lefteqn{ v_b(z)=
\frac{\sin\alpha (\nu_1 +\nu_2 +z)}{\sin\alpha (\nu_2 +z)}\,
\frac{\sin 2\alpha (\nu_2 +z)}{\sin (2\alpha z)} } & & \nonumber \\
& &\times
\frac{\sin\alpha (\nu_1^\prime +\nu_2^\prime +\gamma +z)}
     {\sin\alpha (\nu_2^\prime +\gamma +z)}\,
\frac{\sin 2\alpha (\nu_2^\prime +\gamma +z)}{\sin 2\alpha (\gamma +z)}
\label{vbreduced},
\end{eqnarray}
and
\begin{eqnarray}
\lefteqn{ d_b^+ (z)=
\frac{ (e^{i(\nu_2 +z)};e^{-\beta})_\infty }
     { (e^{i(\nu_1+\nu_2 +z)};e^{-\beta})_\infty }
\frac{ (e^{2iz};e^{-2\beta})_\infty }
     { (e^{2i(\nu_2 +z)};e^{-2\beta})_\infty } } & & \nonumber \\
& &\times
\frac{ (e^{i(\nu_2^\prime +\gamma +z)};e^{-\beta})_\infty }
     { (e^{i(\nu_1^\prime +\nu_2^\prime +\gamma +z)};e^{-\beta})_\infty }
\frac{ (e^{2i(\gamma +z)};e^{-2\beta})_\infty }
     { (e^{2i(\nu_2^\prime +\gamma +z)};e^{-2\beta})_\infty } .
\end{eqnarray}
For the following parameters
$\Delta (x)$ \for{densityb} reduces to
Macdonald's weight function $\Delta_{(R,S)}$
with $R=B_n$, $C_n$ or $BC_n$ and $S=B_n$ or $C_n$:


\setlength{\unitlength}{0.050mm}
\begin{displaymath}
\begin{array}{|c|ccc|} \hline
\begin{picture}(150,100)(35,0)
 \put (5,-20)  {\makebox(75,50)[br]{$S$}}
 \put (0,100){\line(3,-2){215.0}}
 \put (115,30){\makebox(75,50)[b]{$R$}}
\end{picture} & \raisebox{0.3ex}{$B_n$}
              & \raisebox{0.3ex}{$C_n$}
              & \raisebox{0.3ex}{$BC_n$}  \\ [0.5ex] \hline
    &  &\multicolumn{1}{|l|}{}  &    \\ [-1.0ex]
B_n & \nu_1^\prime =\nu_2^\prime =\nu_2 =0
    & \multicolumn{1}{|l|}{\nu_1^\prime =\nu_2^\prime =\nu_1 =0}
    & \nu_1^\prime =\nu_2^\prime =0 \\ [1.0ex] \cline{2-4}
    & &\multicolumn{1}{|l|}{} & \\ [-1.5ex]
\raisebox{2.0ex}{$C_n$}
    & \raisebox{2.0ex}{$ \begin{array}{l}
                          \nu_1^\prime =\nu_1 \\
	                  \nu_2^\prime =\nu_2 =0
                         \end{array} $}
    & \multicolumn{1}{|l|}{
      \raisebox{2.0ex}{$ \begin{array}{l}
                          \nu_1^\prime =\nu_1 =0\\
	                  \nu_2^\prime =\nu_2
                          \end{array} $}  }
    & \raisebox{2.0ex}{$ \begin{array}{l}
                          \nu_1^\prime =\nu_1 \\
	                  \nu_2^\prime =\nu_2
                         \end{array} $} \\ \hline
\end{array}
\end{displaymath}

\begin{center}
{\bf Table.} Special cases associated with admissible pairs $(R,S)$.
\end{center}

\noindent The relation with the parameters employed in Ref. \cite{mac1} reads:
\begin{equation}
t_{\pm e_j \pm e_k}=e^{i\mu},\;\;\;\;\;\;\;\;\;
t_{\pm e_j}= e^{i\nu_1}, \;\;\;\;\;\;\;\;\;
t_{\pm 2e_j}=e^{2i\nu_2}
\end{equation}
and
\begin{equation}
q=\left\{ \begin{array}{ll}
           e^{-\beta } & {\rm if}\;\; S=B_n\; {\rm or}\; S=R=C_n \\
           e^{-\beta /2} & {\rm if}\;\; S=C_n\; {\rm and}\; R=B(C)_n.
         \end{array} \right.
\end{equation}
(In order to verify that
for the above parameters, $\Delta$ \for{densityb} indeed coincides with
the weight functions introduced by
Macdonald, it may be helpful to compare our expressions with Eqs.
(3.1)-(3.5) of \cite{koo3}, since the latter are rather explicit).

Next, we consider the Macdonald polynomials associated with $\Delta_{(R,S)}$.
We distinguish two cases:

{\em i}. $R=(B)C_n$\\
In this case the torus ${\Bbb T}_{R}$($={\Bbb R}^n /(2\pi {\Bbb Z} R^\vee )$)
coincides with ${\Bbb T}$ \for{torus} and
the algebra ${\cal A}^W_{R}$ of W-invariant polynomials on
${\Bbb T}_{R}$ coincides with ${\cal A}^W$.
The fundamental weights of $R=(B)C_n$ read
\begin{equation}
\omega_k =e_1+\cdots +e_k, \;\;\;\;\;\;\;\;\;\;\; k=1,\ldots ,n \label{fundw}
\end{equation}
(our convention regarding the choice of the positive roots agrees
with \cite{bou}).
The cone of dominant weights ${\cal P}^+_{R}$, which consists of
the non-negative
integral combinations of $\omega_k$, $k=1,\ldots ,n$,  coincides
with the cone ${\cal P}^+$ \for{domcone}.
By specializing the parameters as in column 2 and 3 of the above table,
$\{ p_\lambda \}_{\lambda \in {\cal P}^+}$ reduces to an orthogonal
basis of $L^2_W({\Bbb T}_R, \Delta_{(R,S)})$.
Combined with the structure of the expansion \for{kp1}, it now follows that
this basis coincides with
the Macdonald basis $\{ p_{(R,S),\lambda } \}_{\lambda \in {\cal P}^+_R}$.

{\em ii.} $R=B_n$\\
This case is a bit more complicated because
${\Bbb T}$ \for{torus} is merely a subgroup with index two of
the Macdonald torus
${\Bbb T}_{B_n}$($\equiv {\Bbb R}^n /(2\pi {\Bbb Z} B_n^\vee )
                  ={\Bbb R}^n /(2\pi {\Bbb Z} C_n )$).
For $k<n$, the fundamental weights $\omega_k$ of $B_n$ are the same as in
\for{fundw}, but $\omega_n$ is now given by the spin weight
\begin{equation}
\omega_n = (e_1+\cdots +e_n)/2.\label{spinw}
\end{equation}
The cone of dominant weight vectors can be written as
\begin{equation}
{\cal P}^+_{B_n} = \{ \; \lambda +\delta \omega_n \; |\;
\lambda \in {\cal P}^+,\; \delta =0,1\; \} , \label{domBn}
\end{equation}
and the
algebra of W-invariant polynomials on ${\Bbb T}_{B_n}$ is spanned
by the associated monomials:
\begin{equation}
{\cal A}^W_{B_n}= {\rm span} \{ m_{\lambda^\prime}\; |\;
\lambda^\prime \in {\cal P}^+_{B_n}\; \} .
\end{equation}
For $\delta =0$ ($1$)
the function $m_{\lambda +\delta\omega_n}(x)$ \for{monomial} is periodic
(anti-periodic) in $x_j$:
\begin{equation}
x_j \ra x_j +2\pi \; \Longrightarrow \; m_{\lambda +\delta\omega_n}(x)\ra
(-1)^\delta m_{\lambda +\delta \omega_n}(x).\label{anti}
\end{equation}
Combining \for{anti} with the fact that the functions are even in $x_j$ entails
that $m_{\lambda +\omega_n}(x)$ is zero on the hyperplanes
$x_j =\pi$ (mod $2\pi$). Therefore, $m_{\lambda +\omega_n}$ is divisible by
\begin{equation}
m_{\omega_n}(x) =2^n \prod_{1\leq j\leq n} \cos (x_j/2). \label{antigst}
\end{equation}
Thus, we have the following decomposition of ${\cal A}^W_{B_n}$ in periodic
and anti-periodic functions:
\begin{equation}
{\cal A}^W_{B_n} = {\cal A}^W \oplus m_{\omega_n}{\cal A}^W \label{ortdec}
\end{equation}
(with ${\cal A}^W$ as before).
This decomposition is orthogonal with respect to the inner product
on $L^2({\Bbb T}_{B_n}, \Delta _{B_n}dx )$
because $\Delta_{B_n}(x)$ is periodic in $x_j$ with period $2\pi$.

The situation is now as follows: just as the basis
$\{ m_{\lambda^\prime} \}_{\lambda^\prime \in {\cal P}^+_{B_n}}$
the Macdonald basis
$\{ p_{\lambda^\prime ,(B_n ,S)} \}_{\lambda^\prime \in {\cal P}^+_{B_n}}$
of ${\cal A}^W_{B_n}$ splits in periodic and anti-periodic functions
in ${\cal A}^W$ and $m_{\omega_n}{\cal A}^W$, respectively.
By specializing to the first column
of the table, we recover the $B_n$-type Macdonald polynomials that are in the
subspace ${\cal A}^W$ (cf. case {\em i.}, above).
The $B_n$-type polynomials in $m_{\omega_n} {\cal A}^W$
can also be expressed in terms of Koornwinder's polynomials.
To see this, notice that $m_{\omega_n}^2 \Delta_{(B_n,S)}$
coincides with a
weight function of the type \for{densityb}-\for{db+} with parameters
\begin{equation}
\mu_0=\nu_1,\;\;\; \mu_1=i\beta ,\;\;\;
\mu_0^\prime =\left\{ \begin{array}{lcl}
                         0, &\; & S=B_n \\
                         \nu_1, &\; & S=C_n
					  \end{array} \right. , \;\;\;
\mu_1^\prime =0.\label{parcon}
\end{equation}
By multiplying $m_{\omega_n}$ and
the Koornwinder polynomials with parameters as in \for{parcon},
we obtain an orthogonal basis of
$m_{\omega_n}{\cal A}^W$; the latter polynomials coincide with the
anti-periodic Macdonald polynomials.
To be more explicit, we have:
\begin{eqnarray}
p_{\lambda ,(B_n ,S)} &=& p_{\lambda}\;\;\;\;\;\;\;\;\;
\left\{ \begin{array}{llc}
\mu_0 =\nu_1,               &\mu_0^\prime =\mu_1=\mu_1^\prime =0, &
\;\;\;\;\;\;\;  S=B_n \\
\mu_0 =\mu_0^\prime= \nu_1, &\mu_1=\mu_1^\prime =0, & \;\;\;\;\;\;\; S=C_n
       \end{array},\right. \label{pheel} \\
p_{\lambda +\omega_n, (B_n,S)} &=& m_{\omega_n}\: p_{\lambda} \;\;\;
\left\{ \begin{array}{lllc}
\mu_0 =\nu_1,               &\mu_1 =i\beta , &\mu_0^\prime =\mu_1^\prime =0,
& S=B_n \\
\mu_0 =\mu_0^\prime= \nu_1, &\mu_1 =i\beta ,  &\mu_1^\prime =0,
& S=C_n
\end{array} \right. \label{phalf}
\end{eqnarray}
(with $\lambda \in {\cal P}^+$).

\vspace{1ex}
Let us now turn to the corresponding difference operators.
Let $\hat{D}_{1, (R,S)},\ldots ,\hat{D}_{n, (R,S)}$
denote the operators $\hat{D}_1,\ldots ,\hat{D}_n$ \for{Dr1}
with parameters given by the table.
We claim that the polynomials $ p_{\lambda^\prime ,(R,S)}$,
${\lambda^\prime \in {\cal P}^+_R}$,
are joint eigenfunctions of
the operators $\hat{D}_{1, (R,S)},\ldots ,\hat{D}_{n, (R,S)}$.
For $R=(B)C_n$, and for the polynomials \for{pheel}, this is an immediate
consequence of Theorem~\ref{efunc}. For the polynomials \for{phalf} this is
seen
as follows. By conjugating $\hat{D}_{1,(B_n,S)}$ with $m_{\omega_n}$
\for{antigst}, one obtains (up to an additive constant) the operator
$\hat{D}_1$
\for{D1} with parameters \for{parcon}:
\begin{eqnarray}
\lefteqn{m_{\omega_n}^{-1}\: \hat{D}_{1,(B_n,S)}\: m_{\omega_n} } & &
\label{relm} \\
&=&\sum_{\stackrel{1\leq j\leq n}{\eps =\pm 1}} v_{b}(\eps x_j)\prod_{k\neq j}
v_a (\eps x_j+x_k) v_a (\eps x_j-x_k) \left( \frac{\cos \alpha (i\beta +\eps
x_j)}{\cos \alpha\eps x_j} e^{-\eps\beta \hat{\theta}_j}-1 \right) \;\; (\alpha
=1/2) \nonumber \\ &=&
\sum_{\stackrel{1\leq j\leq n}{\eps =\pm 1}}
\frac{\cos\alpha (i\beta +\eps x_j)}{\cos \alpha\eps x_j} v_{b}(\eps x_j)
\prod_{k\neq j} v_a (\eps x_j+x_k) v_a (\eps x_j-x_k)
\left( e^{-\eps\beta \hat{\theta}_j}-1 \right)  \; + \;  const \nonumber
\end{eqnarray}
(calculate the residues to verify the second equality).
Therefore, the polynomials \for{phalf} are eigenfunctions of
$\hat{D}_{1,(B_n,S)}$.
One generalizes this to $\hat{D}_{r, (B_n,S)}$, $r>1$,
via similar reasoning as in the proof of Proposition~\ref{mcApol}:
first one shows that
$\hat{D}_{r, (B_n,S)}$
leaves invariant the space of anti-periodic polynomials
$m_{\omega_n}{\cal A}^W$ (by calculating the residues);
from the asymptotics for Im~$x\ra \infty$ in the positive Weyl chamber,
it then follows that the operator is triangular.
One uses Eq.~\for{relm}, the monotony of
the spectrum of $\hat{D}_{1, (B_n,S)}$
(Lemma~\ref{nondeg}), and the commutativity of the operators to
conclude that the polynomials \for{phalf} are joint eigenfunctions of
$\hat{D}_{1, (B_n,S)},\ldots ,\hat{D}_{n, (B_n,S)}$.

The value of the additive constant in the r.h.s. of \for{relm} can be easily
obtained by comparing the spectrum of the operators on both sides of the
equation (cf. eqs. \for{ev}-\for{roo}, for $r=1$, and \for{parresred}):
\begin{equation}
const = 2\sum _{1\leq j\leq n}\left(
\ch\beta (\rho_j + 1/2)-\ch\beta \rho_j \right)
\end{equation}
with $\rho_j= (n-j)g + k_1 /2$ ($S=B_n$) or
$\rho_j= (n-j)g + k_1$ ($S=C_n$).
In principle one could generalize \for{relm} to an expression that
relates $m_{\omega_n}^{-1}\hat{D}_{r, (B_n,S)}m_{\omega_n}$
to the operators $\hat{D}_1,\ldots ,\hat{D}_r$ with parameters \for{parcon}.
The precise form of these relations can be obtained by
comparing the spectrum of the operators. More precisely, one has to
express $E_{r,n}(\theta )$ \for{ern} in terms of
$E_{1,n}(\theta ),\ldots ,E_{r,n}(\theta )$ with $\rho$ replaced by
$\rho +\omega_n$ and use Theorem~\ref{hc}.

\vspace{1ex}
\noindent\remarks {\em i.} The operator
$\hat{D}_{1, (R,S)}$ coincides (up to a multiplicative constant) with
the Macdonald difference operator
$D_{\pi}$ that is associated with the first fundamental weight
$\omega_1=e_1$ of $S^\vee$. For technical reasons, Macdonald works with a
dilated
root system $S\ra a S$ and the weight lattice of $S^\vee$ is scaled
correspondingly: ${\cal P}_{S^\vee}\ra a^{-1} {\cal P}_{S^\vee}$.
This has as consequence that in comparing with \cite{mac1}
one must multiply the weights of $S^\vee$ with a factor 2 if
$S=C_n$ and $R=B(C)_n$.
Specifically, $\hat{D}_{1,(R,S)} \sim D_{e_1}$ if $S=B_n$ or
$S=R=C_n$, and
$\hat{D}_{1,(R,S)} \sim D_{2e_1}$ if $S=C_n$ and
$R=B(C)_n$.

{\em ii.} The operators $\hat{D}_{1, (R,S)},\ldots , \hat{D}_{n,(R,S)}$
generate
an abelian algebra ${\Bbb D}$ consisting of difference operators that are
simultaneously diagonalized by
$\{ p_{\lambda^\prime ,(R,S)}\}_{\lambda^\prime \in {\cal P}^+_R}$.
Now let $S=C_n$ and consider the operator
\begin{equation}
\hat{D}_{n, (R,C_n)}^\prime \equiv
\sum_{\eps_1,\ldots ,\eps_n =\pm 1}\;
\prod_{1\leq j\leq n} v_b^\prime (\eps_j x_j)\:
\prod_{1\leq j<k\leq n} v_a(\eps_j x_j+ \eps_k x_k)\;
e^{-\beta  (\eps_1 \hat{\theta}_1+\cdots +\eps_n \hat{\theta}_n)/2}\label{Eomn}
\end{equation}
with
\begin{equation}
v_b^\prime (z) \equiv
\frac{ \sin\alpha (\nu_1 +\nu_2 +z)}{\sin \alpha (\nu_2 +z)}\,
\frac{ \sin 2\alpha (\nu_2 +z)}{\sin (2\alpha z)}.
\end{equation}
This \ado{} coincides up to a multiplicative constant with the Macdonald
operator $E_{\pi}$ that is associated with the $n$th fundamental weight
$\omega_n$ of $C^\vee =B_n$; ($\omega_n$ is given by \for{spinw}).
One has (cf. Remark~{\em i})
$\hat{D}_{n, (R,C_n)}^\prime \sim E_{\omega_n}$ if $R=C_n$, and
$\hat{D}_{n, (R,C_n)}^\prime \sim E_{2\omega_n}$ if $R=B(C)_n$.

It follows from \cite{mac1} that
$\hat{D}^\prime_{n,(R,C_n)}$ is diagonalized by
$\{ p_{\lambda^\prime ,(R,C_n)}\}_{\lambda^\prime \in {\cal P}^+_R}$.
Notice however, that $\hat{D}_{n, (R,C_n)}^\prime $ is not in ${\Bbb D}$ but
its
square is (it is easy to see this by examining the eigenvalues and using
Theorem~\ref{hc}).

{\em iii.} From a group-theoretic perspective, Theorem~\ref{hc}
amounts to saying that ${\Bbb D}$ is isomorphic to ${\Bbb R}[{\cal P}]^W$ (the
W-invariant part of the (real) group algebra over the lattice
${\cal P}={\Bbb Z}^n$).
If $S=B_n$, then one has ${\cal P}_{S^\vee}={\cal P}$; so
${\Bbb D}_{(R,B_n)}\equiv {\Bbb D}\cong {\Bbb R}[{\cal P}_{B_n^\vee}]^W$.
If $S=C_n$, then ${\cal P}$ is a subgroup of ${\cal P}_{S^\vee}$ with index
two;
so ${\Bbb D}$ is isomorphic to a subalgebra of
${\Bbb R}[{\cal P}_{C_n^\vee}]^W$.
In the latter case, one can extend ${\Bbb D}$ to an algebra
${\Bbb D}_{(R,C_n)}$
that is isomorphic
to ${\Bbb R}[{\cal P}_{C_n^\vee}]^W$ by replacing the generator
$\hat{D}_{n, (R,C_n)}$ by $\hat{D}^\prime_{n, (R,C_n)}$ \for{Eomn}.

\subsection{The Root System $D_n$}
We conclude by briefly sketching the state of affairs for $R=D_n$.
(The interested reader should not have much difficulty to supply
missing proofs by comparing with the previous
subsection).
Put
\begin{equation}
\mu_0=\mu_1=\mu_0^\prime=\mu_1^\prime=0. \label{Dnpar}
\end{equation}
Then $d_b(z)=1$ and $\Delta (x)$ \for{densityb} reduces to Macdonald's
$D_n$-type weight function. (Again the correspondence of parameters is via
Eq. \for{psubsta}). Also $v_b(z)=1$ and $\hat{D}_1$ reduces to Macdonald's
operator
$D_{\omega_1}$ associated with the first fundamental weight $\omega_1=e_1$.

For $R=D_n$, the torus ${\Bbb T}_R={\Bbb R}^n / (2\pi {\Bbb Z}R^\vee )$ is
the same as for $R=B_n$. The Weyl group, however, is smaller: only
an even number of sign flips of the variables $x_j$, $j=1,\ldots , n$, is
allowed. For $k=1,\ldots ,n-2$, the fundamental weights $\omega_k$ of $D_n$
are the same as in \for{fundw}, but $\omega_{n-1}$ and $\omega_n$
are now given by the half-spin weights
\begin{equation}
\omega_{n-1} =  (e_1+\cdots +e_{n-1}-e_n) /2,\;\;\;\; \;\;\; \;\;\;\;
\omega_{n} = (e_1+\cdots +e_{n-1}+e_n)/2.\label{halfspw}
\end{equation}
It is not hard to see that the cone of dominant weights ${\cal P}_{D_n}^+$
generated
by $\omega_k$, $k=1,\ldots ,n$, consists of the vectors
\begin{equation}
(\lambda +\delta \omega_n)_\eps \equiv
(\lambda_1 +\delta /2,\ldots ,\lambda _{n-1}+\delta /2,\,
\eps (\lambda _{n}+\delta /2) )
\end{equation}
with
\begin{equation}
\lambda \in {\cal P}^+,\;\;\;\;\;
\delta =0,1,\;\;\;\;\;\;
\eps = \pm 1.
\end{equation}
The Macdonald polynomials $p_{\lambda^\prime ,D_n}$,
$\lambda^\prime \in {\cal P}^+_{D_n}$ form an orthogonal basis
of $L^2_{W_{D_n}}({\Bbb T}_{D_n},\Delta_{D_n} dx)$.
By combining the polynomials associated with $(\lambda +\delta\omega_n)_+$
and $(\lambda +\delta\omega_n)_-$ one obtains polynomials that are even in
$x_j$, $j=1,\ldots ,n$. These are related to Koornwinder's polynomials
in the following way (cf. Eqs. \for{pheel}, \for{phalf}):
\begin{equation}
\begin{array}{rllc}
p_{\lambda_+ ,D_n}=p_{\lambda_- ,D_n}=& p_{\lambda},&
\mu_0=\mu_1=\mu_0^\prime =\mu_1^\prime =0, &
{\rm if}\; \lambda_n =0, \\
p_{\lambda_+ ,D_n}+p_{\lambda_- ,D_n}=& p_{\lambda}, &
\mu_0=\mu_1=\mu_0^\prime =\mu_1^\prime =0, &
{\rm if}\; \lambda_n >0,  \\
p_{(\lambda +\omega_n)_+ ,D_n}+ p_{(\lambda +\omega_n)_- ,D_n}=
& m_{\omega_n}\; p_{\lambda}, &
\mu_0=\mu_0^\prime =\mu_1^\prime =0,\; \mu_1 =i\beta .&
\end{array} \label{Dpols}
\end{equation}
In the second and the third line of the above formula one obtains a sum
of $D_n$ polynomials rather than the polynomials themselves.
Nevertheless, Eq.~\for{Dpols} determines
the $D_n$ polynomials uniquely. This is because
flipping the sign of one of the $x_j$'s in
$p_{(\lambda +\delta \omega_n)_+ ,D_n}$ results in
$p_{(\lambda +\delta \omega_n)_- ,D_n}$. Consequently,
the coefficients of the expansion of
$p_{(\lambda +\delta \omega_n)_\eps ,D_n}$ in $D_n$-type monomial
symmetric functions are determined in terms of the coefficients occurring
in \for{kp1}.

As regards the difference operators with parameters \for{Dnpar},
the algebra ${\Bbb D}$ consists of commuting \ados{}
with the $D_n$-type polynomials as joint eigenfunctions. One can extend
${\Bbb D}$ to an algebra that is isomorphic to
${\Bbb R}[{\cal P}_{D_n^{\vee }}]^{W_{D_n}}$ by replacing the generators
$\hat{D}_{n-1}$ and $\hat{D}_{n}$ by
\begin{eqnarray}
\hat{D}_{n-1}^\prime &=&
\sum_{\stackrel{ \eps_1,\ldots ,\eps_n =\pm 1}{\eps_1\cdots
\eps_n =-1}} \;\;\;\prod_{1\leq j< k\leq n} v_a(\eps_j x_j +\eps_k x_k)\;
e^{-\beta (\eps_1 \hat{\theta}_1+\cdots +\eps_n \hat{\theta}_n)/2}, \\
\hat{D}_{n}^\prime &=&
\sum_{\stackrel{ \eps_1,\ldots ,\eps_n =\pm 1}{\eps_1\cdots
\eps_n =+1}} \;\;\;\prod_{1\leq j< k\leq n} v_a(\eps_j x_j +\eps_k x_k)\;
e^{-\beta  (\eps_1 \hat{\theta}_1+\cdots + \eps_n \hat{\theta}_n)/2}.
\end{eqnarray}
These operators are proportional to Macdonald's operators $E_{\omega_{n-1}}$
and
$E_{\omega_n}$, which are associated with the half-spin weights \for{halfspw}.

\renewcommand{\theequation}{A.\arabic{equation}}
\renewcommand{\thesection}{A}
\markright{Appendix A}
\addtocontents{toc}{\protect\vspace{-2ex}}
\addcontentsline{toc}{section}{Appendix A: Cancellation of Poles}
\setcounter{equation}{0}
\setcounter{defn}{0}
\section*{Appendix A: Cancellation of Poles}
In this appendix we prove two results, which were needed to demonstrate
that $\hat{D}_r$ maps ${\cal A}^W$ into itself.
It was claimed in the proof of Proposition~\ref{inv} (Section~\ref{sectria})
that the following expression:
\begin{equation}
V^1_J\, V^3_{J;J^c}
\sum_{\stackrel{\emptyset\subsetneqq J_1\subsetneqq \cdots \subsetneqq J_s=J}
               {1\leq s\leq |J|}}
(-1)^s\prod_{1\leq s^\prime\leq s}V^2_{J_{s^\prime}\setminus J_{s^\prime -1}}
V^3_{J_{s^\prime}\setminus J_{s^\prime -1};J\setminus J_{s^\prime}}
\; [m_\lambda (x+2\gamma e_{J_1})-m_\lambda (x)]\label{regeq}
\end{equation}
(with $\lambda \in {\cal P}^+$ and $J\subset \{ 1,\ldots ,n\}$,
$J_0\equiv \emptyset$), is regular
both at $x_1=-\gamma$ (pole of type $I$) and at $x_1=-x_j-2\gamma$,
$j=2,\ldots ,n$ (poles of type $II$).
The terms of \for{regeq} have poles due to zeros in the
denominators of the coefficients (cf. \for{VJ1}-\for{VJ3} and
\for{va},\for{vb}).
Recall that we assume that the parameters
$\gamma$, $\mu$, $\mu_\delta$, $\mu^\prime_\delta$ ($\delta =0,1$)
and the variables $x_2,\ldots ,x_n$ are chosen
in such a way that these poles are simple.
In the next two lemmas we prove the above regularity claims,
thereby completing the proof of Proposition~\ref{inv}.

Before turning to the details, let us outline the idea of the proof.
Equation \for{regeq} consists of a sum of terms of the type
\begin{equation}
(-1)^s\, V^1_J\, V^3_{J;J^c}
\prod_{1\leq s^\prime\leq s}\, V^2_{B_{s^\prime}}\,
V^3_{B_{s^\prime};J\setminus J_{s^\prime}}\;
[m_\lambda (x+2\gamma e_{B_1})-m_\lambda (x)]\label{terms}
\end{equation}
where the index sets  $B_{s^\prime}\subset J$, $s^\prime =1,\ldots ,s$ denote
the blocks of
the partition of the cell $J$:
\begin{equation}
B_{s^\prime}\equiv J_{s^\prime} \setminus J_{s^\prime -1},\;\;\;\;\;\;\;
1\leq s^\prime \leq s.
\end{equation}
The terms \for{terms} are associated  with the sequences
\begin{equation}
\emptyset \subsetneqq J_1\subsetneqq J_2\subsetneqq \cdots \subsetneqq J_s=J,
\;\;\;\;\;\;\;\;\;\; (1\leq s\leq |J|) \label{seq}
\end{equation}
(with the cell $J$ fixed).
Each term in \for{regeq} corresponds to a sequence \for{seq}.
We will construct an involutive operation $\sigma$ ($\sigma^2 =$id)
on the collection of sequences \for{seq} in such a way that the terms
associated with a sequence and its image under $\sigma$ have opposite
residue. Therefore, the poles in \for{regeq} cancel in pairs.

\begin{lem}[pole of type $I$]\label{pole1}\hfill

\noindent Let $\gamma$, $\mu$, $\mu_\delta$, $\mu^\prime_\delta$ ($\delta
=0,1$)
and $x_2,\ldots ,x_n$ be such
that the terms in \for{regeq} have only simple poles.
Then \for{regeq} is regular as a function of $x_1$ at $x_1=-\gamma$.
\end{lem}
\begin{prf}
First note that the lemma is trivial if $1\notin J$, because in that case
$V^1_J$ does not depend on $x_1$. But if $1\in J$, then
$V_J^1$ gives rise to a pole at $x_1=-\gamma$ in \for{regeq}.

Assume $1\in J$ and let $B_{s_1}$ denote the block of the cell $J$ that
contains the index 1.
We define the following map $\sigma$ on the collection of
sequences \for{seq}:
\begin{itemize}
\item[1.A] If $|B_{s_1}|>1$, then $\sigma$ maps
\for{seq} to the sequence
\begin{equation}
\emptyset \subsetneqq J_1 \subsetneqq J_2\subsetneqq \cdots\subsetneqq
J_{s_1-1}\subsetneqq \; J_{s_1}\setminus \{ 1\} \; \subsetneqq
J_{s_1}\subsetneqq
\cdots \subsetneqq  J_s=J. \label{imseq}
\end{equation}
\item[1.B] If $|B_{s_1}|=1$ and $s_1 >1$, then $\sigma$ maps
\for{seq} to the sequence
\begin{equation}
\emptyset \subsetneqq J_1 \subsetneqq J_2\subsetneqq \cdots \subsetneqq
J_{s_1-2}\subsetneqq  J_{s_1}\subsetneqq
\cdots \subsetneqq  J_s=J.
\end{equation}
\item[2.] If $|B_{s_1}|=1$ and $s_1 =1$ (i.e.
$J_1 =\{ 1\}$), then $\sigma$ maps sequence \for{seq} to itself.
\end{itemize}
Phrased in words: unless $B_{s_1}=\{ 1\}$, the map $\sigma$ pulls the index $1$
out of
$B_{s_1}$ and places it in a newly created block,
which is sandwiched between $B_{s_1}\setminus \{ 1\}$ and
$B_{s_1+1}$  (case 1.A); when $B_{s_1}$ contains only the index $1$,
then $\sigma$ merges the blocks
$B_{s_1}=\{ 1\} $ and $B_{s_1-1}$ if $s_1>1$ (case 1.B) or,
if $s_1=1$, then it leaves the sequence \for{seq} unchanged (case 2.).

Thus defined, $\sigma$ is indeed an involution on the collection
of sequences \for{seq}: the cases 1.A and 1.B are inverse to each other
(see Fig. 1. below).

\begin{center}

\setlength{\unitlength}{0.8mm}
\begin{picture}(140,50)(5,0)


\put(0,37){\thicklines\line(1,0){50}}
\put(0,30){\makebox(15,7){$\cdots\cdots$}}
\put(14,30){\framebox(22,7){$j_1,\cdots j_p,\; 1$}}
\put(35,30){\makebox(15,7){$\cdots\cdots$}}
\put(0,30){\thicklines\line(1,0){50}}

\put(70,37){\makebox(10,7){1.A}}
\put(60,37){\vector(1,0){30}}
\put(90,30){\vector(-1,0){30}}
\put(70,23){\makebox(10,7){1.B}}

\put(100,37){\thicklines\line(1,0){50}}
\put(100,30){\makebox(15,7){$\cdots\cdots$}}
\put(114,30){\framebox(22,7){$j_1,\cdots j_p\;\;\: 1$}}
\put(135,30){\makebox(15,7){$\cdots\cdots$}}
\put(100,30){\thicklines\line(1,0){50}}
\put(131,30){\line(0,1){7}}


\put(5,10){\makebox(45,7){$\cdots\cdots\cdots\cdots\cdots\cdots\cdots$}}
\put(0,10){\framebox(5,7){1}}
\put(0,17){\thicklines\line(1,0){50}}
\put(0,10){\thicklines\line(1,0){50}}

\put(70,13.5){\vector(-1,0){10}}
\put(70,10){\makebox(10,7){2.}}
\put(80,13.5){\vector(1,0){10}}

\put(105,10){\makebox(45,7){$\cdots\cdots\cdots\cdots\cdots\cdots\cdots$}}
\put(100,10){\framebox(5,7){1}}
\put(100,17){\thicklines\line(1,0){50}}
\put(100,10){\thicklines\line(1,0){50}}

\end{picture}

\end{center}
\begin{center}
{\bf Fig. 1.} A graphical representation of the map $\sigma$.
\end{center}

We claim that in the first situation (i.e. 1.A or 1.B) the pole at
$x_1=-\gamma$
in the term \for{terms} (which is associated with \for{seq}) cancels against
the
pole in the term corresponding with the  $\sigma$-image of the sequence
\for{seq}. To see this, we may assume that we are in situation 1.A.
One
obtains the term corresponding to the sequence \for{imseq} from \for{terms} by
making the substitutions
\begin{eqnarray} s&\ra & s+1 ,\\
V^2_{B_{s_1}}&\ra &
V^2_{B_{s_1}\setminus \{ 1\} }\: V^2_{\{ 1\}} =
V^2_{B_{s_1}\setminus \{ 1\} } ,\\
V^3_{B_{s_1};J\setminus J_{s_1}}&\ra &
V^3_{B_{s_1}\setminus \{ 1\} ; (J\setminus J_{s_1})\cup \{ 1\} }\:
V^3_{\{ 1\} ; J\setminus J_{s_1}} \nonumber \\
& &= V^3_{B_{s_1};J\setminus J_{s_1}}\:
V^3_{B_{s_1}\setminus \{ 1\} ; \{ 1\} } ,\\
m_\lambda (x+2\gamma e_{B_1})&\ra &
m_\lambda (x+2\gamma e_{B_1}-2\gamma \delta_{1,s_1} e_1 )
\end{eqnarray}
( $\delta _{j,k}$ denotes the Kronecker symbol).
This substitution in \for{terms} amounts to replacing the part
\begin{eqnarray}
\lefteqn{\left( \left. V^2_{B_{s_1}} \right/
V^2_{B_{s_1}\setminus \{ 1\} } \right) \:
[m_\lambda (x+2\gamma e_{B_1})-m_\lambda (x)]=} & & \label{ter1} \\
& &\prod_{j\in B_{s_1}\setminus \{ 1\} }
v_a(x_j+x_1)\, v_a(x_j+x_1+2\gamma )\:
[m_\lambda (x+2\gamma e_{B_1})-m_\lambda (x)] \nonumber
\end{eqnarray}
by
\begin{eqnarray}
\lefteqn{-V^3_{B_{s_1}\setminus \{ 1\} ; \{ 1\} }\:
[m_\lambda (x+2\gamma e_{B_1}-2\gamma\delta_{1,s_1}e_1)
-m_\lambda (x)]=} & & \label{ter2} \\
& &-\prod_{j\in B_{s_1}\setminus \{ 1\} }
v_a(x_j+x_1)\, v_a(x_j-x_1)\:
[m_\lambda (x+2\gamma e_{B_1}-2\gamma\delta_{1,s_1}e_1)-m_\lambda (x)].
\nonumber
\end{eqnarray}
At $x_1=-\gamma$ the r.h.s. of \for{ter1} and \for{ter2}
differ only by sign.
Hence, the residues at $x_1=-\gamma$ cancel.

If we are in situation 2., i.e. $B_1=\{ 1\}$, then \for{terms} is regular at
$x_1=-\gamma$ because the pole in $V^1_J$ is compensated by a zero in
the difference of the two monomial symmetric functions:
\begin{equation}
\left[ m_\lambda (x+2\gamma e_1)-m_\lambda (x)\right]_{(x_1 =-\gamma )} =0.
\end{equation}

This shows that the total residue at $x_1=-\gamma$ of the sum \for{regeq}
vanishes, which completes the proof of the lemma.
\end{prf}

\begin{lem}[poles of type $II$]\label{pole2}\hfill

\noindent Let $\gamma$, $\mu$, $\mu_\delta$, $\mu^\prime_\delta$ ($\delta
=0,1$)
and $x_2,\ldots ,x_n$ be such
that the terms in \for{regeq} have only simple poles.
Then \for{regeq} is regular as a function of $x_1$ at
$x_1=-x_j-2\gamma$, $j=2,\ldots , n$.
\end{lem}
\begin{prf}
The proof is very similar to that of Lemma~\ref{pole1}. Fix a
$j\in \{ 2,\ldots , n\}$. The lemma is trivial if $J$ does not contain the pair
$\{ 1,j\}$, because in that case all terms of \for{regeq} are regular at
$x_1=-x_j-2\gamma$.

Assume for the remaining part of the proof that $\{ 1,j\} \subset J$. Let
$\sigma_j$ be the following map on the collection of sequences \for{seq}:
\begin{itemize}
\item[1.A]  If the pair $\{ 1,j\}$ is contained in one of the blocks,
say $B_{s_{j}}$, of sequence \for{seq}, and
$|B_{s_{j}}|>2$, then $\sigma_j$ maps sequence \for{seq}
to
\begin{equation}
\emptyset \subsetneqq J_1\subsetneqq \cdots \subsetneqq
J_{s_{j}-1}\subsetneqq\; J_{s_{j}}\setminus \{ 1,j\} \;\subsetneqq J_{s_{j}}
\subsetneqq \cdots \subsetneqq J_s=J.
\end{equation}
\item[1.B] If one of the blocks of sequence \for{seq},
say $B_{s_{j}}$, equals $\{ 1,j\}$, and $s_{j} >1$,
then $\sigma_j$ maps sequence \for{seq} to
\begin{equation}
\emptyset \subsetneqq J_1\subsetneqq \cdots \subsetneqq J_{s_{j}-2}\subsetneqq
J_{s_{j}} \subsetneqq \cdots \subsetneqq J_s=J.
\end{equation}
\item[2.] If $B_1 =\{ 1,j\}$, then $\sigma_j$ maps sequence \for{seq} to
itself.
\item[3.] If the pair $\{ 1,j\}$ is not contained in any of the blocks
of \for{seq}, then $\sigma_j$ maps the sequence \for{seq} to itself.
\end{itemize}
It is clear that $\sigma_j$ is an involution, the cases 1.A and 1.B
are inverse to each other (see Fig 2. below).
\begin{center}

\setlength{\unitlength}{0.8mm}

\begin{picture}(140,75)(5,0)


\put(0,57){\thicklines\line(1,0){50}}
\put(0,50){\makebox(15,7){$\cdots\cdot$}}
\put(12,50){\framebox(26,7){$j_1,\cdots j_p,\; 1,j$}}
\put(35,50){\makebox(15,7){$\cdots\cdot$}}
\put(0,50){\thicklines\line(1,0){50}}

\put(70,57){\makebox(10,7){1.A}}
\put(60,57){\vector(1,0){30}}
\put(90,50){\vector(-1,0){30}}
\put(70,43){\makebox(10,7){1.B}}

\put(100,57){\thicklines\line(1,0){50}}
\put(100,50){\makebox(15,7){$\cdots\cdot$}}
\put(112,50){\framebox(26,7){$j_1,\cdots j_p\;\;\: 1,j$}}
\put(135,50){\makebox(15,7){$\cdots\cdot$}}
\put(100,50){\thicklines\line(1,0){50}}
\put(129.5,50){\line(0,1){7}}


\put(5,30){\makebox(45,7){$\cdots\cdots\cdots\cdots\cdots\cdots\cdots$}}
\put(0,30){\framebox(7,7){$1,j$}}
\put(0,37){\thicklines\line(1,0){50}}
\put(0,30){\thicklines\line(1,0){50}}

\put(70,33.5){\vector(-1,0){10}}
\put(70,30){\makebox(10,7){2.}}
\put(80,33.5){\vector(1,0){10}}

\put(105,30){\makebox(45,7){$\cdots\cdots\cdots\cdots\cdots\cdots\cdots$}}
\put(100,30){\framebox(7,7){$1,j$}}
\put(100,37){\thicklines\line(1,0){50}}
\put(100,30){\thicklines\line(1,0){50}}


\put(0,17){\thicklines\line(1,0){50}}
\put(0,10){\makebox(5,7){$\cdot\cdot$}}
\put(5,10){\framebox(17,7){$i_1\cdot\cdot i_p, 1$}}
\put(22,10){\makebox(5,7){$\cdot\cdot$}}
\put(27,10){\framebox(17,7){$j_1\cdot\cdot j_{p^\prime}, j$}}
\put(45,10){\makebox(5,7){$\cdot\cdot$}}
\put(0,10){\thicklines\line(1,0){50}}

\put(70,13.5){\vector(-1,0){10}}
\put(70,10){\makebox(10,7){3.}}
\put(80,13.5){\vector(1,0){10}}

\put(100,17){\thicklines\line(1,0){50}}
\put(100,10){\makebox(5,7){$\cdot\cdot$}}
\put(105,10){\framebox(17,7){$i_1\cdot\cdot i_p, 1$}}
\put(122,10){\makebox(5,7){$\cdot\cdot$}}
\put(127,10){\framebox(17,7){$j_1\cdot\cdot j_{p^\prime}, j$}}
\put(145,10){\makebox(5,7){$\cdot\cdot$}}

\put(100,10){\thicklines\line(1,0){50}}

\end{picture}

\end{center}
\vspace{2ex}
\begin{center}
{\bf Fig. 2.} A graphical representation of the map $\sigma_j$.
\end{center}

Consider situation 1., assuming case 1.A. The application of
$\sigma_j$ boils down to making the following
substitutions in the associated term \for{terms}:
\begin{eqnarray}
s&\ra & s+1 \\
V^2_{B_{s_{j}}} &\ra&
V^2_{B_{s_{j}}\setminus \{ 1,j\} }\:
V^2_{\{ 1,j\} } \\
V^3_{B_{s_{j}};\, J\setminus J_{s_{j}}}&\ra&
V^3_{B_{s_{j}};\, J\setminus J_{s_{j}}}\:
V^3_{B_{s_{j}};\, \{ 1,j\} } .
\end{eqnarray}
These substitutions amount to the following change in the term \for{terms}:
replace the part
\begin{equation}
\prod_{\stackrel{ k\in B_{s_{j}}}{k\neq 1,j}}
v_a(x_k+x_1)\, v_a(x_k+x_1+2\gamma )\, v_a(x_k+x_j)\, v_a(x_k+x_j+2\gamma
)\label{lpart} \end{equation}
by
\begin{equation}
-\prod_{\stackrel{ k\in B_{s_{j}}}{k\neq 1,j}}\label{rpart}
v_a(x_k+x_1)\, v_a(x_k-x_1)\, v_a(x_k+x_j)\, v_a(x_k-x_j).
\end{equation}
At $x_1=-x_j-2\gamma$, \for{lpart} and \for{rpart} differ only by sign.
Consequently, the residues at  $x_1=-x_j-2\gamma$ of the corresponding terms in
\for{regeq} add up to zero.

In situation 2. the pole in the coefficient of \for{terms}, which is caused by
$V^2_{B_{1}}$, cancels against the zero in the difference
of the monomial symmetric functions:
\begin{equation}
\left[ m_\lambda (x+2\gamma e_{\{ 1,j\} })-m_\lambda (x)
\right]_{(x_1=-x_j-2\gamma )}=0.
\end{equation}

In situation 3. the denominator of the coefficient of \for{terms} has no
zero at $x_1=-x_j-2\gamma$, so the term is regular at
$x_1=-x_j-2\gamma$.

We conclude from the above analysis that the total residue at
$x_1=-x_j-2\gamma$
in the sum \for{regeq} is zero, thus completing the proof of the lemma.
\end{prf}

\renewcommand{\theequation}{B.\arabic{equation}}
\renewcommand{\thesection}{B}
\markright{Appendix B}
\addtocontents{toc}{\protect\vspace{-2ex}}
\addcontentsline{toc}{section}{Appendix B: Two Combinatorial Lemmas}
\setcounter{equation}{0}
\setcounter{defn}{0}
\section*{Appendix B: Two Combinatorial Lemmas}
In this appendix we prove two technical results used in the main text,
which have a bearing on the eigenvalues of our difference operators.
In Lemma~\ref{reclem1} we solve a certain linear system;
its solution resulted in explicit formulas for the eigenvalues of
$\hat{D}_r$ (Proposition~\ref{Deigv}).
Lemma~\ref{reclem2} deals with a recursion relation, which
helped us in
obtaining explicit information concerning the behavior of the
eigenvalues for $\beta\ra 0$ (Proposition~\ref{evbto0}).

\begin{lem}\label{reclem1}
The functions
\begin{equation}
F_{m ,p}= (-1)^{p}
\sum_{1\leq i_1\leq \cdots \leq i_{p}\leq m-p+1}
t_{i_1}t_{i_2}\cdots t_{i_{p}},
\;\;\;\;\;\;\;\;\;\;\; 1\leq p \leq m  ,\label{solu}
\end{equation}
form the unique solution of the linear system
\begin{equation}
\sum_{\stackrel{J\subset \{ 1,\ldots , n\} ,\, |J|=s}{0\leq s\leq r}}
\left( \prod_{j\in J} t_j \right) F_{n-s,r-s}=0,
\;\;\;\;\;\;\;\;\;\; 1\leq r \leq n  ,\label{linsyst}
\end{equation}
with the convention
\begin{equation}
F_{m ,0}=1,\;\;\;\;\;\;\;\;\;\;   m =0,1,2,\ldots.\label{inc}
\end{equation}
\end{lem}
\begin{prf}
After splitting off the term in \for{linsyst} corresponding to $s=0$
and bringing all other terms to the r.h.s. of the equation, one arrives
at a recursion relation for $F_{n,r}$:
\begin{equation}
F_{n,r}=- \sum_{\stackrel{J\subset \{ 1,\ldots , n\} ,\; |J|=s}{1\leq s\leq r}}
\left( \prod_{j\in J} t_j \right) F_{n-s,r-s},
\;\;\;\;\;\;\;\;\;\; 1\leq r \leq n.\label{receq}
\end{equation}
It is clear that \for{receq} with condition \for{inc} determines $F_{n,r}$
uniquely (use induction on $r$).
Hence, the system \for{linsyst}, \for{inc} has a unique solution.

In order to prove that this solution is indeed given by Eq. \for{solu},
we must show that the expression
\begin{equation}
\sum_{\stackrel{J\subset \{ 1,\ldots , n\} ,\; |J|=s}{0\leq s\leq r}}
(-1)^s \left( \prod_{j\in J} t_j \right)
\left(\sum_{1\leq i_1\leq \cdots \leq i_{r-s}\leq n-r+1}
t_{i_1}t_{i_2}\cdots t_{i_{r-s}}\right) \label{van}
\end{equation}
vanishes identically (for $1\leq r\leq n$).
To this end we observe that
\for{van} consists of a sum of monomials in the variables
$t_1,\ldots ,t_n$ of the type
\begin{equation}\label{monterm}
(-1)^s\; (t_{j_1}t_{j_2}\cdots t_{j_s})\times
(t_{i_1}t_{i_2}\cdots t_{i_{r-s}}),
\end{equation}
which correspond to pairs of the form
\begin{equation}
\{ (j_1,j_2,\ldots , j_s),\; (i_1,i_2,\ldots , i_{r-s}) \}
,\;\;\;\;\;\;\;\;\;\;
0\leq s\leq r,\label{tup1}
\end{equation}
subject to the condition
\begin{equation}
1\leq j_1< j_2 < \cdots < j_s\leq n,\;\;\;\;\;
1\leq i_1\leq i_2 \leq \cdots \leq i_{r-s}\leq n-r+1.\label{tup2}
\end{equation}
We shall now show that these monomials cancel in pairs.

Let $\sigma$ be the following operation
defined on the above collection of
pairs \for{tup1} with condition \for{tup2}:
\begin{itemize}
\item[A.] If $i_1<j_1$ or $s=0$, then
\begin{equation}
\{ (j_1,j_2,\ldots , j_s),\; (i_1,i_2,\ldots , i_{r-s}) \}
\stackrel{\sigma}{\longrightarrow}
\{ (i_1,j_1,\ldots , j_s),\; (i_2,\ldots , i_{r-s}) \} ,
\end{equation}
\item[B.] If $i_1\geq j_1$ or $s=r$, then
\begin{equation}
\{ (j_1,j_2,\ldots , j_s),\; (i_1,i_2,\ldots , i_{r-s}) \}
\stackrel{\sigma}{\longrightarrow}
\{ (j_2,\ldots , j_s),\; (j_1,i_1,\ldots , i_{r-s}) \} .
\end{equation}
\end{itemize}
Roughly speaking, $\sigma$ compares the first entries of the two elements
constituting the pair \for{tup1} and moves the smallest of these two
to the first entry of the other element.
One easily verifies that:
first, $\sigma$ is well defined in the sense that the image of \for{tup1}
is again a pair satisfying \for{tup2};
second, $\sigma$ is an involution ($\sigma^2 =id$), the cases {A.}
and {B.} being inverse to each other.

For the associated monomial \for{monterm}, acting with $\sigma$ amounts
to an increase (case {A.}) or a decrease (case {B.}) of the number $s$
by one, i.e. it flips the sign of the corresponding monomial.
Therefore, combining the term \for{monterm} associated with a pair \for{tup1}
with the one associated with its image under $\sigma$ entails
the vanishing of the sum \for{van}, which completes the proof.
\end{prf}

\vspace{1ex}
\noindent\remark If one replaces the upper bound $n-r+1$ of the second
summation
in \for{van} by $n$, then, for $r=n$, expression \for{van}
also vanishes. (Indeed, the above proof again applies).
In this case the vanishing of \for{van} amounts to a well-known
relation between the elementary
symmetric functions and the complete symmetric functions
(see e.g. \cite{mac0}).

\begin{lem}\label{reclem2}
The function
\begin{eqnarray}
\lefteqn{E_{r,n}(t_1,\ldots , t_n;p_r,\ldots , p_n)=} & &\label{solE} \\
& & \sum_{0\leq s\leq r}(-1)^{r+s}
\left( \sum_{1\leq j_1<\cdots < j_s\leq n} t_{j_1}\cdots t_{j_s} \right)
\left(\sum_{r\leq i_1\leq\cdots \leq i_{r-s}\leq n} p_{i_1}\cdots
p_{i_{r-s}}\right), \;\;\;\;\;\;\;\; 1\leq r\leq n ,\nonumber
\end{eqnarray}
is the unique solution of the recursion relation
\begin{eqnarray}
\begin{array}{r}
E_{r,n}(t_1,\ldots , t_n;p_r,\ldots , p_n)=
(t_n-p_n)E_{r-1,n-1}(t_1,\ldots , t_{n-1};p_r,\ldots , p_n)+  \\
E_{r,n-1}(t_1,\ldots , t_{n-1};p_r,\ldots , p_{n-1}),
\end{array} \nonumber \\
\hfill  1\leq r\leq n \label{rrel}
\end{eqnarray}
with the convention
\begin{equation}
E_{0,n}\equiv 1, \;\;\;\;\;\;\;\;\;\;
E_{r,n}\equiv 0 \;\; {\rm if}\; n<r. \label{intc}
\end{equation}
\end{lem}
\begin{prf}
It is clear that \for{rrel} with condition \for{intc} determines $E_{r,n}$
uniquely (use induction on $n$).
After splitting up the sum in \for{solE} in three parts,
it becomes apparent that \for{solE} indeed solves Eq. \for{rrel}:

\begin{tabular}{llr}
{\em i}. &terms with $j_s=n$:     &$t_n E_{r-1,n-1}(t_1,\ldots , t_{n-1};
                                                          p_r,\ldots , p_n)$;
\\
{\em ii}. &terms with $j_s<n$ and
                         $i_{r-s}=n$: &$-p_n E_{r-1,n-1}(t_1,\ldots , t_{n-1};
                                                          p_r,\ldots , p_n)$;
\\

{\em iii}. &terms with $j_s<n$ and
                         $i_{r-s}<n$: &$E_{r,n-1}(t_1,\ldots , t_{n-1};
                                                       p_r,\ldots , p_{n-1})$.
\end{tabular}
\end{prf}

\vspace{1ex}
\noindent\remark In some cases Lemma~\ref{reclem2} can be used to obtain
alternative
expressions for $E_{r,n}$. For instance, one easily verifies
with the aid of relation \for{rrel} that if
\begin{equation}
p_r,p_{r+1},\ldots ,p_n=\overline{p},
\end{equation}
then
\begin{equation}
E_{r,n} =
\sum_{\stackrel{J\subset \{ 1,\ldots , n\}}{|J|=r}}
\left( \prod_{j\in J} ( t_j- \overline{p} ) \right) =
S_r(t_1-\overline{p},\ldots ,t_n-\overline{p} ) .
\end{equation}
In particular,
$E_{n,n}(t_1,\ldots ,t_n;p_n)=(t_1-p_n)\, (t_2-p_n)\cdots (t_n-p_n)$.

\renewcommand{\theequation}{C.\arabic{equation}}
\renewcommand{\thesection}{C}
\markright{Appendix C}
\addtocontents{toc}{\protect\vspace{-2ex}}
\addcontentsline{toc}{section}{Appendix C: $\hat{D}({\cal
A}^W)=0\Longrightarrow
\hat{D}=0$}
\setcounter{equation}{0}
\setcounter{defn}{0}
\section*{Appendix C: $\hat{D}({\cal A}^W)=0 \Longrightarrow \hat{D}=0$}
In this appendix we present a result due to
S. N. M. Ruijsenaars \cite{rui4}. It shows
that if an \ado{} or \pdo{} is zero on all symmetric functions
in ${\cal A}^W$, then its coefficients must be zero.
This fact was used in Section~\ref{diacom}
to show that the operators $\hat{D}_1,\ldots ,\hat{D}_n$
commute (Theorem~\ref{commuting}), and again in Section~\ref{oplim}
to conclude that for $\beta \ra 0$ one obtains the $BC_n$-type hypergeometric
\pdos{} of Heckman and Opdam.

Let
\begin{equation}
\hat{T}_{\kappa}\equiv {\hat{t}_1}^{\kappa_1}\cdots {\hat{t}_n}^{\kappa_n},
\;\;\;\;\;\;\;\;\;\;\;
{\hat{t}_j}\equiv
\left\{ \begin{array}{ccc}
         e^{-\beta\hat{\theta}_j} & (\mbox{\ado{}}) & (\beta >0) \\
		 \hat{\theta}_j            & (\mbox{\pdo{}}) &
		\end{array} \right. ,
\end{equation}
with $\kappa =(\kappa_1,\ldots ,\kappa_n)$ in ${\Bbb R}^n$ or ${\Bbb N}^n$ in
the
\ado{} or \pdo{} case, respectively.
The \ados{}/\pdos{} of interest are of the form:
\begin{equation}
\hat{D}=\sum_{1\leq r\leq M} V_r(x)\: \hat{T}_{\kappa^{(r)}},\label{operator}
\end{equation}
with the $n$-dimensional vectors
$\kappa^{(1)},\ldots ,\kappa^{(M)}$ distinct, and the coefficient
functions
\begin{equation}
V_r:{\cal U}\subset {\Bbb R}^n \longrightarrow {\Bbb C},\;\;\;\;\;\;\;\;\;\;
1\leq r\leq M,
\end{equation}
continuous on an open dense set ${\cal U}\subset {\Bbb R}^n$.

\begin{prp}[$\hat{D}({\cal A}^W)=0 \Rightarrow \hat{D}=0$]\label{vanish}\hfill

\noindent Let $\hat{D}$ be an \ado{}/\pdo{} of the form
\for{operator}. If
\begin{equation}
\hat{D}\, m_\lambda  =0, \;\;\;\;\;\;\;\;\;\; \forall \lambda \in {\cal P}^+,
\end{equation}
then
\begin{equation}
V_1(x)=V_2(x)=\cdots =V_M(x)=0,\;\;\;\;\;\;\; \forall x \in {\cal U}.
\end{equation}
\end{prp}
{\em Proof} (\cite{rui4})\newline
Introduce the following vector-valued functions:
\begin{eqnarray}
t_\lambda : {\cal U}\ra {\Bbb C}^M,& &\;\;\;\;\;\;
x\stackrel{t_\lambda}{\longmapsto}
\left(
(\hat{T}_{\kappa^{(1)}}\, m_{\lambda})(x),\ldots ,(\hat{T}_{\kappa^{(M)}}\,
m_{\lambda})(x)
\right) ,\;\;\; \lambda\in {\cal P}^+, \\
v : {\cal U}\ra {\Bbb C}^M,& &\;\;\;\;\;\;
x\stackrel{v}{\longmapsto}
\left( \overline{V_1(x)},\ldots ,\overline{V_M(x)} \right) .\label{vvector}
\end{eqnarray}
The fact that $m_\lambda$ is in the kernel of $\hat{D}$ translates itself
geometrically in the orthogonality of $t_\lambda$ and $v$:
\begin{equation}
\hat{D}\, m_\lambda =0 \Longleftrightarrow
(t_\lambda ,v)=0.
\end{equation}
We will assume $v\not\equiv 0$ and derive a contradiction.
Let $\lambda^{(1)},\ldots , \lambda^{(M)}$ be vectors in ${\cal P}^+$. One has
$v\perp t_{\lambda^{(1)}},\ldots , t_{\lambda^{(M)}}$.
Therefore, the vectors $t_{\lambda^{(1)}}(x),\ldots , t_{\lambda^{(M)}}(x)$
must be
linearly dependent for all $x\in {\cal U}$ for which $v(x)\neq 0$.
Since $v(x)$ is continuous in $x$, there exists an open ball
${\cal B}\subset {\cal U}$ on which $v(x)\neq 0$. The fact that the
vectors $t_{\lambda^{(s)}}(x)$, $1\leq s\leq M$ are real-analytic in $x$
then entails
\begin{equation}
{\rm det}\left( {\rm Col}[ t_{\lambda^{(1)}}(x),\ldots , t_{\lambda^{(M)}}(x)]
\right) = 0,\;\;\;\;\;\;\;\; \forall x\in {\cal U}.\label{detnul}
\end{equation}
We will now show that an appropriate choice of the vectors
$\lambda^{(1)},\ldots , \lambda^{(M)}$ contradicts the vanishing of the above
determinant.

Let $\lambda \in {\cal P}^+$ and $y\in {\Bbb R}^n$ be such that
\begin{equation}
\lambda_1\geq\lambda_2\geq \cdots \geq \lambda_n >0,
\;\;\;\;\;\; \;\;\; \;\;\;\;\;\;\;
y_1>y_2>\cdots > y_n>0. \label{pstl}
\end{equation}
{}From the asymptotics (cf. Eq. \for{as1})
\begin{equation}
(\hat{T}_{\kappa}\, m_\lambda )(x)|_{x=iRy} \sim
 \tau_{\kappa ,\lambda}\: e^{ (\lambda ,y)R}, \;\;\;\;\;\;\;\;\;\;\;\;\;
R\ra \infty ,
\end{equation}
with
\begin{equation}
\tau_{\kappa ,\lambda }=\left\{
\begin{array}{lc}
e^{\beta (\kappa ,\lambda )} & (\mbox{\ado{}}) \\
(-1)^{|\kappa |}\: (\lambda_1)^{\kappa_1}\cdots (\lambda_n)^{\kappa_n} &
(\mbox{\pdo{}})
\end{array} \right. ,\label{taud}
\end{equation}
one derives:
\begin{equation}
\lim_{R\ra \infty}\;
e^{i (\lambda ,x)}\: t_{\lambda}(x) |_{x=iRy}\; =\;
(\tau_{\kappa^{(1)},\lambda},\ldots , \tau_{\kappa^{(M)},\lambda}).\label{limt}
\end{equation}

Pick the vector $\lambda$ (subject to condition \for{pstl}) in such a way that
\begin{equation}
\tau_{\kappa^{(r)},\lambda}\neq \tau_{\kappa^{(p)},\lambda},\;\;\;\;\;\;\;\;\;
1\leq r< p \leq M.\label{distinct}
\end{equation}
That such a $\lambda$ exists follows in the \ado{} case from
\for{taud} and the fact that the vectors $\kappa^{(1)},\ldots , \kappa^{(M)}$
are distinct;
in the \pdo{} case one can pick distinct prime numbers for the components of
$\lambda$.

We use $\lambda =(\lambda_1,\ldots , \lambda_n)$ to form the vectors
$\lambda^{(1)},\ldots , \lambda^{(M)}$ in the following way:
\begin{equation}
\lambda^{(s)} =
\left\{ \begin{array}{lc}
(s-1)\, (\lambda_1,\ldots , \lambda_n) & (\mbox{\ado{}}) \\ [0.5ex]
\left( (\lambda_1)^{s-1},\ldots , (\lambda_n)^{s-1}\right) & (\mbox{\pdo{}})
\end{array} \right. ,\;\;\;\;\;\;\;\; 1\leq s\leq M.\label{labs}
\end{equation}
On the one hand, Eqs. \for{detnul} and \for{limt} imply
\begin{equation}
\tau \equiv
\left|
\begin{array}{ccc}
\tau_{\kappa^{(1)},\lambda^{(1)}} &\cdots & \tau_{\kappa^{(1)},\lambda^{(M)}}
\\
\vdots                    &\ddots & \vdots                    \\
\tau_{\kappa^{(M)},\lambda^{(1)}} &\cdots & \tau_{\kappa^{(M)},\lambda^{(M)}}
\end{array}\right| =0. \label{vandet}
\end{equation}
On the other hand, for the above choice of the vectors
$\lambda^{(s)}$ \for{labs} , $\tau$ is a Vandermonde determinant:
\begin{equation}
\tau_{\kappa^{(r)},\lambda^{(s)}}=(\tau_{\kappa^{(r)},\lambda})^{s-1},
\;\;\;\;\;\;\;\;\;\;\;
1\leq s\leq M.
\end{equation}
Therefore,
\begin{equation}
\tau =\prod_{1\leq r< p\leq M}
(\tau_{\kappa^{(p)},\lambda}-\tau_{\kappa^{(r)},\lambda} ).\label{vdm}
\end{equation}
Because $\lambda$ is chosen such that
$\tau_{\kappa^{(r)},\lambda}\neq \tau_{\kappa^{(p)},\lambda}$
if $r\neq p$,
it follows from Eq.~\for{vdm} that the determinant $\tau \neq 0$,
contradicting \for{vandet}.

Hence, $v$ \for{vvector} must be zero.
\epr

{\noindent\bf Acknowledgments.} The author would like to thank
S. N. M. Ruijsenaars for many helpful conversations and useful suggestions
during the whole process that led to this paper.
Thanks are also due to T. H. Koornwinder for explanation of his results.

\addtocontents{toc}{\protect\vspace{-2ex}}
\addcontentsline{toc}{section}{References}
\markright{References}

\end{document}